\documentclass[12pt]{article}
\usepackage{amsmath,amsfonts,amssymb,amsthm}
\begin{document}

\newcommand{\1}{{{\bf 1}}}
\newcommand{\id}{{\rm id}}
\newcommand{\Hom}{{\rm Hom}\,}
\newcommand{\End}{{\rm End}\,}
\newcommand{\Res}{{\rm Res}\,}
\newcommand{\Image}{{\rm Im}\,}
\newcommand{\Ind}{{\rm Ind}\,}
\newcommand{\Aut}{{\rm Aut}\,}
\newcommand{\Ker}{{\rm Ker}\,}
\newcommand{\gr}{{\rm gr}}
\newcommand{\Der}{{\rm Der}\,}

\newcommand{\Z}{\mathbb{Z}}
\newcommand{\Q}{\mathbb{Q}}
\newcommand{\C}{\mathbb{C}}
\newcommand{\N}{\mathbb{N}}

\newcommand{\g}{\mathfrak{g}}
\newcommand{\h}{\mathfrak{h}}
\newcommand{\wt}{{\rm wt}\;}
\newcommand{\CR}{\mathcal{R}}
\newcommand{\D}{\mathcal{D}}
\newcommand{\E}{\mathcal{E}}
\newcommand{\Lie}{\mathcal{L}}
\newcommand{\z}{\bf{z}}
\newcommand{\bflam}{\bf{\lambda}}

\def \<{\langle} 
\def \>{\rangle}
\def \be{\begin{equation}\label}
\def \ee{\end{equation}}
\def \bex{\begin{exa}\label}
\def \eex{\end{exa}}
\def \bl{\begin{lem}\label}
\def \el{\end{lem}}
\def \bt{\begin{thm}\label}
\def \et{\end{thm}}
\def \bp{\begin{prop}\label}
\def \ep{\end{prop}}
\def \br{\begin{rem}\label}
\def \er{\end{rem}}
\def \bc{\begin{coro}\label}
\def \ec{\end{coro}}
\def \bd{\begin{de}\label}
\def \ed{\end{de}}

\newtheorem{thm}{Theorem}[section]
\newtheorem{prop}[thm]{Proposition}
\newtheorem{coro}[thm]{Corollary}
\newtheorem{conj}[thm]{Conjecture}
\newtheorem{exa}[thm]{Example}
\newtheorem{lem}[thm]{Lemma}
\newtheorem{rem}[thm]{Remark}
\newtheorem{de}[thm]{Definition}
\newtheorem{hy}[thm]{Hypothesis}
\makeatletter
\@addtoreset{equation}{section}
\def\theequation{\thesection.\arabic{equation}}
\makeatother
\makeatletter

\begin{center}{\Large \bf On certain categories of modules
for affine Lie algebras}\\
\vspace{0.5cm}
Haisheng Li\footnote{Partially supported by a NSA grant
and a grant from Rutgers Research Council}\\
Department of Mathematical Sciences, Rutgers University, Camden, NJ 08102\\
and\\
Department of Mathematics, Harbin Normal University, Harbin, China
\end{center}

\begin{abstract}
In this paper, we re-examine certain integrable modules of
Chari-Presslely for an (untwisted) affine Lie algebra $\hat{\g}$
by exploiting basic formal variable techniques.
We define and study two categories ${\mathcal{E}}$ and ${\mathcal{C}}$ of 
$\hat{\g}$-modules using generating functions,
where ${\mathcal{E}}$ contains evaluation modules and
${\mathcal{C}}$ unifies highest weight modules, evaluation
modules and their tensor product modules, and we classify
integrable irreducible $\hat{\g}$-modules in categories
${\mathcal{E}}$ and ${\mathcal{C}}$.
\end{abstract}

\baselineskip=16pt

\section{Introduction}
Let $\g$ be a finite-dimensional simple 
Lie algebra equipped with the Killing form
$\<\cdot,\cdot\>$ which is suitably normalized.
Associated to the pair $(\g,\<\cdot,\cdot\>)$ 
we have the (untwisted) affine Lie
algebra $\hat{\g}$ (without the degree derivation added). 
For affine algebras, a very important class of modules 
is the class of highest weight modules (cf. [K1])
in the well known category ${\mathcal{O}}$, where
highest weight integrable (irreducible) modules
(of nonnegative integral levels) have been the main focus.
We also have another class of modules,
called evaluation modules (of level zero)
associated with a finite number of
$\g$-modules and with the same number of nonzero complex numbers,
studied by Chari and Presslely in [CP2] (cf. [CP1], [CP3]). 
Furthermore, Chari and Presslely in [CP2] studied the first time
the tensor product module of an integrable highest 
weight $\hat{\g}$-module with
a (finite-dimensional) evaluation $\hat{\g}$-module associated with
finite-dimensional irreducible $\g$-modules and distinct 
nonzero complex numbers.
(Such a tensor product module is integrable
as the tensor product module of any two integrable modules
is integrable.) 
A surprising result, proved in [CP2], is that 
such a tensor product module is also irreducible. 
In this way, a new family of integrable
$\hat{\g}$-modules were constructed. 

We know that integrable highest 
weight $\hat{\g}$-modules are exactly the irreducible integrable modules
in the well-known category ${\mathcal{O}}$ (see [K1]) and 
that irreducible integrable evaluation modules are exactly 
the finite-dimensional irreducible modules (see [C], [CP2]).
In view of this, naturally one would want to
find a canonical characterization for the new integrable modules,
{\em instead of presenting them as tensor product modules}.
This is part of our motivation for this paper.
Part of our motivation is to look for canonical connections 
of the new integrable modules with modules and fusion rules
for affine vertex operator algebras.

In this paper, we give a canonical characterization of 
the new integrable modules using generating functions and formal calculus. 
Notice that highest weight $\hat{\g}$-modules belong to a bigger class of 
modules called restricted modules, where
a $\hat{\g}$-module $W$ is said to be {\em restricted} (cf. [K1]) if
for any $a\in \g,\; w\in W$, 
$(a\otimes t^{n})w=0$ for $n$ sufficiently large.
In terms of generating functions, a $\hat{\g}$-module $W$ is restricted 
if and only if $a(x)w\in W((x))$ for $a\in \g,\; w\in W$, 
where $a(x)=\sum_{n\in \Z}(a\otimes t^{n})x^{-n-1}$ (the generating function).
For an evaluation module $U$ (see [CP2]), we show that
there is a nonzero polynomial $p(x)$ such that
$p(x)a(x)u=0$ for $a\in \g,\; u\in U$. Then $p(x)a(x)v\in (W\otimes U)((x))$ for
$a\in \g,\; v\in W\otimes U$. Motivated by these facts,
we define a category $\E$ to consist of
$\hat{\g}$-modules $W$ such that there exists a nonzero polynomial $p(x)$ 
such that $p(x)a(x)w=0$ for $a\in \g,\; w\in W$ and we define
a category ${\mathcal{C}}$ to consist of
$\hat{\g}$-modules $W$ such that there exists a nonzero polynomial $f(x)$ 
such that $f(x)a(x)w\in W((x))$ for $a\in \g,\; w\in W$.
Then category $\E$ contains all the evaluation modules and 
category ${\mathcal{C}}$ contains all the restricted modules,
the evaluation modules and their tensor products, so that 
category ${\mathcal{C}}$ unifies all the mentioned modules.
In this paper we prove that the irreducible integrable $\hat{\g}$-modules
in the category $\E$ are exactly the finite-dimensional irreducible
evaluation modules up to isomorphism. (This result is analogous and closely related
to a result of Chari-Presslely [C], [CP2].)
It was proved in [DLM] that every restricted integrable $\hat{\g}$-module
is a direct sum of highest weight irreducible integrable $\hat{\g}$-modules.
As our main result of this paper we prove that
the irreducible integrable $\hat{\g}$-modules
in the category ${\mathcal{C}}$ up to isomorphism are exactly the 
tensor product modules of highest weight irreducible integrable
$\hat{\g}$-modules with finite-dimensional irreducible evaluation modules.
The key to our main result is a factorization result which states that
every irreducible representation of $\hat{\g}$ 
in the category ${\mathcal{C}}$ 
can be factorized canonically as the product of two representations
of $\hat{\g}$ such that the first representation 
defines a restricted module and the second one 
defines a module in the category $\E$. The proof of this factorization
uses formal calculus in an essential way.

It is well known (cf. [Li2], [LL]) that restricted $\hat{\g}$-modules 
are closely related to affine vertex operator algebras and their modules.
But the tensor product $\hat{\g}$-modules in the category ${\mathcal{C}}$
is not a module for the affine vertex operator algebra.
In this paper, by using a result of \cite{li-duke} we show that
if $W$ and $W_{1}$ are highest weight integrable
irreducible $\hat{\g}$-modules of the same level and $U(z)$ is a
finite-dimensional evaluation module, 
$\Hom_{\hat{\g}}(W\otimes U(z),\overline{W_{2}})$
gives the fusion rule of a certain type as generally defined in [FHL] 
in terms of vertex operator algebras and their modules.

In this paper, most of the results are proved in 
the generality that $\g$ is only assumed to be of countable
dimension, so those results in fact hold for toroidal Lie algebras.

This paper is organized as follows: In Section 2, in the first half we review 
the definitions and examples of restricted modules and evaluation modules 
for affine Lie algebras and we recall certain results of Chari-Presslely.
In the second half we define categories $\E$ and ${\mathcal{C}}$ and
we give slight generalizations of Chari-Presslely's results.
In Section 3, we classify the irreducible integrable modules
in the categories $\E$ and ${\mathcal{C}}$. In Section 4, 
we give a connection between the tensor product module of a
highest weight irreducible integrable module with an evaluation module 
and fusion rules of certain types.

\section{Categories ${\mathcal{R}}$, ${\mathcal{E}}$ and
${\mathcal{C}}$ of modules for affine Lie algebras}

In this section we review the definitions and examples of
restricted modules and evaluation modules 
for an affine Lie algebra $\hat{\g}$. We define
a category ${\mathcal{E}}$ of $\hat{\g}$-modules, including
evaluation modules, and we define
a category ${\mathcal{C}}$ of $\hat{\g}$-modules, including
restricted modules, evaluation modules and their tensor product modules.
We give a generalization of certain results of 
Chari and Presslely ([C], [CP2]) with
a different proof using formal calculus.

First let us fix some formal variable notations
(see [FLM], [FHL], [LL]).
Throughout this paper, 
$x, x_{1}, x_{2},\dots$
are independent mutually commuting formal variables. 
We shall typically use $z, z_{1},z_{2},\dots$ for complex numbers.
For a vector space $U$, $U[[x_{1}^{\pm 1},\dots,x_{n}^{\pm 1}]]$
denotes the space of all formal (possibly doubly infinite) series
in $x_{1},\dots,x_{n}$ with coefficients in $U$, $U((x_{1},\dots,x_{n}))$
denotes the space of all formal (lower truncated) Laurent series 
in $x_{1},\dots,x_{n}$ with coefficients in $U$
and $U[[x_{1},\dots,x_{n}]]$ denotes the space of all formal
(nonnegative) powers series 
in $x_{1},\dots,x_{n}$ with coefficients in $U$.

\br{rassociativity}
{\em As it was pointed out in [FLM] (cf. [LL]), in formal calculus, 
associativity law and cancellation law
for products of formal series {\em do not} hold in general, 
but they {\em do hold}
if all the involved (sub)products exist. For example, 
if $a(x)\in U((x))$, $f(x), g(x)\in \C((x))$, we have
$$f(x)(g(x)a(x))=(f(x)g(x))a(x).$$
For $a(x),b(x)\in U((x))$, if
$h(x)a(x)=h(x)b(x)$ for some $h(x)\in \C((x))$,
then $a(x)=b(x)$.}
\er

We shall use the traditional binomial expansion convention: 
For $m\in \Z$,
\begin{eqnarray}
(x_{1}\pm x_{2})^{m}=\sum_{i\ge 0}\binom{m}{i}(\pm 1)^{i}x_{1}^{m-i}x_{2}^{i}
\in \C[x_{1},x_{1}^{-1}][[x_{2}]].
\end{eqnarray}

Recall from [FLM] the formal delta function
\begin{eqnarray}
\delta(x)=\sum_{n\in \Z}x^{n}\in \C[[x,x^{-1}]].
\end{eqnarray}
Its fundamental property is that
\begin{eqnarray}
f(x)\delta(x)=f(1)\delta(x)\;\;\;\mbox{ for }f(x)\in \C[x,x^{-1}].
\end{eqnarray}
For any nonzero complex number $z$,
\begin{eqnarray}
\delta\left(\frac{z}{x}\right)=\sum_{n\in \Z}z^{n}x^{-n}\in \C[[x,x^{-1}]]
\end{eqnarray}
and we have
\begin{eqnarray}
f(x)\delta\left(\frac{z}{x}\right)
=f(z)\delta\left(\frac{z}{x}\right)\;\;\;\mbox{ for }f(x)\in \C[x,x^{-1}].
\end{eqnarray}
In particular,
\begin{eqnarray}
(x-z)\delta\left(\frac{z}{x}\right)=0.
\end{eqnarray}

Let $\g$ be a Lie algebra (not necessarily finite-dimensional) equipped with 
a nondegenerate symmetric invariant bilinear form $\<\cdot,\cdot\>$,
fixed throughout this section. 
Let $\hat{\g}$ be the corresponding (untwisted) affine Lie algebra, i.e.,
\begin{eqnarray}
\hat{\g}=\g\otimes \C[t,t^{-1}]\oplus \C {\bf k}
\end{eqnarray}
with the defining commutator relations
\begin{eqnarray}\label{edef-affine-comp}
[a\otimes t^{m},b\otimes t^{n}]
=[a,b]\otimes t^{m+n}+m\<a,b\>\delta_{m+n,0}{\bf k} 
\;\;\;\mbox{ for }a,b\in\g,\; m,n\in \Z\label{eaffine1}
\end{eqnarray}
and with ${\bf k}$ as a nonzero central element.
A $\hat{\g}$-module $W$ is said to be of {\em level} $\ell$ in $\C$ 
if the central 
element ${\bf k}$ acts on $W$ as the scalar $\ell$.
By the {\em standard untwisted
affine algebra} $\hat{\g}$ we mean 
the affine Lie algebra $\hat{\g}$ with
$\g$ a finite-dimensional simple Lie algebra
and with $\<\cdot,\cdot\>$ the normalized Killing form
so that the squared length of the longest roots is $2$.

For $a\in \g$, form the generating function
\begin{eqnarray}
a(x)=\sum_{n\in \Z}(a\otimes t^{n})x^{-n-1}\in \hat{\g} [[x,x^{-1}]].
\end{eqnarray}
In terms of generating functions the defining relations
(\ref{eaffine1}) exactly amount to
\begin{eqnarray}\label{edef-commutator}
[a(x_{1}),b(x_{2})]
=[a,b](x_{2})x_{2}^{-1}\delta\left(\frac{x_{1}}{x_{2}}\right)+\<a,b\>{\bf k}
\frac{\partial}{\partial x_{2}}x_{2}^{-1}\delta\left(\frac{x_{1}}{x_{2}}\right).
\end{eqnarray}

Following the tradition (cf. [FLM], [LL]), for $a\in \g,\; n\in\Z$ 
we shall use $a(n)$ 
for the corresponding operator associated to $a\otimes t^{n}$
on  $\hat{\g}$-modules.
We now introduce the category ${\mathcal{R}}$ of the so-called 
restricted modules for the affine algebra $\hat{\g}$.
A $\hat{\g}$-module $W$ is said to be {\em restricted} (cf. [K1]) if
for any $w\in W,\; a\in \g$,
\begin{eqnarray}\label{erestricted-module}
a(n)w=0\;\;\;\mbox{ for }n\;\;\mbox{sufficiently large}.
\end{eqnarray}
We define a $\Z$-grading $\hat{\g}=\coprod_{n\in \Z}\hat{\g}_{(n)}$ by
\begin{eqnarray}
\hat{\g}_{(0)}=\g\oplus \C {\bf k}\;\;\mbox{ and }\;\;
\hat{\g}_{(n)}=\g\otimes t^{-n}\;\;\;\mbox{ for }n\ne 0,
\end{eqnarray}
making $\hat{\g}$ a $\Z$-graded Lie algebra.
It is clear that any $\N$-graded $\hat{\g}$-module 
is automatically a restricted module.

Let $U$ be a $\g$-module and let $\ell$ be any complex number.
Let ${\bf k}$ act on $U$ as the scalar $\ell$ and let $\g\otimes t\C[t]$ 
act trivially, making $U$ a $(\g\otimes \C [t]\oplus \C{\bf k})$-module.
Form the following induced $\hat{\g}$-module
\begin{eqnarray}\label{einduced-module}
M_{\hat{\g}}(\ell,U)=U(\hat{\g})\otimes_{U(\g\otimes \C[t]\oplus \C{\bf k})}U.
\end{eqnarray}
Endow $U$ with zero degree, making $M_{\hat{\g}}(\ell,U)$ 
an $\N$-graded $\hat{\g}$-module.
This in particular implies that $M_{\hat{\g}}(\ell,U)$ is 
a restricted $\hat{\g}$-module.
This $\hat{\g}$-module is commonly called the {\em Weyl module}
or the {\em generalized Verma module associated with $\g$-module $U$}.
If $\g$ is a finite-dimensional simple Lie algebra and if $U$ is a
(highest weight) Verma $\g$-module, then $M_{\hat{\g}}(\ell,U)$
is isomorphic to a (highest weight) Verma $\hat{\g}$-module (cf. [K1]).
Furthermore, any (highest weight) Verma $\hat{\g}$-module 
is isomorphic to a module of the form $M_{\hat{\g}}(\ell,U)$.
A homomorphic image of a Verma $\hat{\g}$-module is
called a {\em highest weight module}.
Then the category ${\mathcal{R}}$ contains all the highest weight 
modules for the standard affine Lie algebra $\hat{\g}$.

For the affine Lie algebra $\hat{\g}$, we also have another family of
$\hat{\g}$-modules, called the evaluation modules (see [CP2]).
Let $U$ be a $\g$-module and let $z$ be a nonzero
complex number. Define an action of $\hat{\g}$ on $U$ by
\begin{eqnarray}
a(n)\cdot u&=&z^{n}(au)\;\;\;\mbox{ for }a\in \g,\; n\in \Z,\\
{\bf k}\cdot U&=&0.
\end{eqnarray}
Then $U$ equipped with the defined action is a $\hat{\g}$-module (of
level zero) (see [CP2]), which is denoted by $U(z)$.
If $U$ is an irreducible $\g$-module, it is clear that $U(z)$ is an
irreducible $\hat{\g}$-module.
More generally, let $U_{1},\dots,U_{r}$ be $\g$-modules and let
$z_{1},\dots,z_{r}$ be nonzero complex numbers. Then
$U=U_{1}\otimes \cdots \otimes U_{r}$ is a $\hat{\g}$-module
where ${\bf k}$ acts as zero and  
\begin{eqnarray}\label{etensorproductaction1}
a(n)(u_{1}\otimes \cdots\otimes u_{r})
=\sum_{i=1}^{r}z_{i}^{n}(u_{1}\otimes \cdots\otimes au_{i}\otimes 
\cdots \otimes u_{r})
\end{eqnarray}
for $a\in \g,\; n\in \Z,\; u_{i}\in U_{i}$.
This module is nothing but the tensor product
$\hat{\g}$-module $\otimes_{i=1}^{r}U_{i}(z_{i})$.
Such a $\hat{\g}$-module is called an {\em evaluation module}.
The following results are due to Chari and Presslely (see [C] and [CP2]):

\bt{tfinite-dimensional-module}
Let $\g$ be a finite-dimensional simple Lie algebra.
Let $U_{1},\dots,U_{r}$ 
be (finite-dimensional) irreducible $\g$-modules and let $z_{1},\dots,z_{r}$
be distinct nonzero complex numbers. Then 
$\otimes_{i=1}^{r}U_{i}(z_{i})$ is a (finite-dimensional) irreducible 
$\hat{\g}$-module of level zero. 
Furthermore, every finite-dimensional irreducible $\hat{\g}$-module 
is isomorphic to such a $\hat{\g}$-module.
\et

\br{rintegrable-module}
{\em Let $\g$ be a finite-dimensional simple Lie algebra.
We fix a Cartan subalgebra $\h$ and 
denote by $\Delta$ the set of roots, so that 
$\g=\h\oplus \sum_{\alpha\in \Delta}\g_{\alpha}$. 
We also fix a choice of set  $\Delta_{+}$ of positive roots and
denote by $\theta$ the highest long root.
Let $\<\cdot,\cdot\>$ be the normalized Killing form on $\g$ such that 
$\<\theta,\theta\>=2$. A $\hat{\g}$-module $W$ is said to be
{\em integrable} (see [K1], [K2]) 
if $\g_{\alpha}(n)$ acts locally nilpotently on $W$
for $\alpha\in \Delta,\; n\in \Z$.
Then (see [K2]) the subalgebra
$\h\oplus \C {\bf k}$ of $\hat{\g}$ acts semisimply on $W$. 
If $W$ is an irreducible integrable $\hat{\g}$-module, the central element
${\bf k}$ acts on $W$ as a scalar $\ell$ in $\N$.
A {\em singular vector } of a $\hat{\g}$-module $W$
of level $\ell$ is a (nonzero) $\h$-eigenvector $u$ such that
$a(n)u=0$ for $a\in \g,\; n>0$ and $a(0)u=0$ for $a\in \g_{+}$.
A known fact is that the submodule of an integrable $\hat{\g}$-module
generated by a singular vector is irreducible (see [K1]). }
\er

The following result was established by Chari and Presslely in [CP2]:

\bt{tcp}
Let $\hat{\g}$ be a standard affine Lie algebra
(with $\g$ a finite-dimensional simple Lie algebra and $\<\cdot,\cdot\>$ 
the normalized Killing form).
Let $W$ be an irreducible highest weight integrable $\hat{\g}$-module,
let $U_{1},\dots,U_{r}$ be finite-dimensional 
irreducible $\g$-modules and let 
$z_{1},\dots,z_{r}$ be distinct nonzero complex numbers. 
Then the tensor product $\hat{\g}$-module 
$W\otimes U_{1}(z_{1})\otimes \cdots \otimes U_{r}(z_{r})$ is irreducible.
\et

Note that a restricted $\hat{\g}$-module is defined 
canonically by the property (\ref{erestricted-module}) while
typical evaluation $\hat{\g}$-modules are finite-dimensional. 
Then naturally one would want to find a canonical characterization
for the new family of (tensor product) $\hat{\g}$-modules 
$W\otimes U_{1}(z_{1})\otimes \cdots \otimes U_{r}(z_{r})$.
In the following we give a characterization 
in terms of generating functions.

First consider restricted $\hat{\g}$-modules (in the category
${\mathcal{R}}$).
Note that the condition (\ref{erestricted-module}) amounts to that
\begin{eqnarray}\label{erestricted-condition}
a(x)w\in W((x))\;\;\;\mbox{ for }a\in \g, \; w\in W.
\end{eqnarray}
That is, a $\hat{\g}$-module $W$ is restricted if and only if
\begin{eqnarray}
a(x)\in \Hom (W,W((x)))\;\;\;\mbox{ for }a\in \g.
\end{eqnarray}

Then we consider evaluation $\hat{\g}$-modules.
Let $U_{1},\dots, U_{r}$ be $\g$-modules and $z_{1},\dots,z_{r}$ 
nonzero complex numbers.
For $a\in \g,\; u_{i}\in U_{i}(z_{i})=U_{i}$, 
writing (\ref{etensorproductaction1}) in terms of generating functions, 
we have
\begin{eqnarray}\label{eevaluation-generating}
a(x)(u_{1}\otimes \cdots\otimes u_{r})
&=&\sum_{n\in \Z}\sum_{i=1}^{r}z_{i}^{n}x^{-n-1}
(u_{1}\otimes \cdots \otimes au_{i}\otimes \cdots \otimes u_{r})\nonumber\\
&=&\sum_{i=1}^{r}x^{-1}\delta\left(\frac{z_{i}}{x}\right)
(u_{1}\otimes \cdots \otimes au_{i}\otimes \cdots \otimes u_{r}).
\end{eqnarray}
Since $(x-z_{i})\delta\left(\frac{z_{i}}{x}\right)=0$ for $i=1,\dots,r$, we get
$(x-z_{1})\cdots (x-z_{r})a(x)(u_{1}\otimes \cdots\otimes u_{r})=0$.
In view of this and (\ref{erestricted-condition}) we immediately have:

\bl{lchari1}
Let $U_{1},\dots, U_{r}$ be $\g$-modules and 
let $z_{1},\dots,z_{r}$ be nonzero complex numbers. 
Then on the tensor product
$\hat{\g}$-module $U_{1}(z_{1})\otimes
\cdots \otimes U_{r}(z_{r})$,
\begin{eqnarray}
(x-z_{1})\cdots (x-z_{r})a(x)=0\;\;\;\mbox{ for }a\in \g.
\end{eqnarray}
Furthermore, for any restricted $\hat{\g}$-module $W$, we have
\begin{eqnarray}
(x-z_{1})\cdots (x-z_{r})a(x)\in \Hom (M,M((x)))\;\;\;\mbox{ for }a\in \g,
\end{eqnarray}
where $M$ denotes the tensor product $\hat{\g}$-module
$W\otimes U_{1}(z_{1})\otimes \cdots \otimes U_{r}(z_{r})$.
\el

Motivated by Lemma \ref{lchari1}, we define the following two categories:

\bd{dcategoryE-D}
{\em We define a category ${\mathcal{E}}$
to consist of $\hat{\g}$-modules $W$ for which there exists a
nonzero polynomial $p(x)\in \C[x]$ such that
\begin{eqnarray}
p(x)a(x)w=0\;\;\;\mbox{ for }a\in \g,\; w\in W.
\end{eqnarray}
We define category ${\mathcal{C}}$ 
to consist of 
$\hat{\g}$-modules $W$ such that there exists a nonzero polynomial 
$p(x)$ such that 
\begin{eqnarray}
p(x)a(x)\in \Hom (W,W((x)))\;\;\;\mbox{ for }a\in \g.
\end{eqnarray}}
\ed

\br{runification}
{\em In view of Lemma \ref{lchari1}, 
all the evaluation $\hat{\g}$-modules belong to
the category ${\mathcal{E}}$ and all the restricted
$\hat{\g}$-modules, evaluation $\hat{\g}$-modules and
tensor products of restricted
$\hat{\g}$-modules with evaluation $\hat{\g}$-modules
belong to the category ${\mathcal{C}}$.}
\er

\br{rchari}
{\em In [C], Chari defined a category $\tilde{\mathcal{O}}$ of 
$\hat{\g}$-modules and classified all the irreducible modules
and all the integrable modules in this category.
Furthermore, Chari-Presslely proved in [CP2] that
irreducible integrable modules in category $\tilde{\mathcal{O}}$
are exactly the finite-dimensional evaluation modules up to isomorphism.
The categories $\E$ and $\tilde{\mathcal{O}}$ are closely related, 
but they are different.}
\er

\bl{l-level-E}
The central element ${\bf k}$ acts as zero on any
$\hat{\g}$-module in the category ${\mathcal{E}}$.
\el

\begin{proof} Let $W$ be a $\hat{\g}$-module in the category
${\mathcal{E}}$ with a nonzero polynomial $p(x)$ such that
$p(x)u(x)=0$ on $W$ for $u\in \g$. If $p(x)$ is of degree zero, 
we have $u(x)=0$ for $u\in \g$, i.e., $u(n)=0$ for $u\in \g,\; n\in \Z$.
In view of the commutator relation (\ref{edef-affine-comp}) we see that
${\bf k}$ must be zero on $W$.
Assume that $p(x)$ is not a constant, that is, $p'(x)\ne 0$.
Let $a,b\in \g$ be such that $\<a,b\>=1$.
(Notice that $\<\cdot,\cdot\>$ is assumed to be nondegenerate.)
Using the commutator relations (\ref{edef-commutator}) we get
\begin{eqnarray}
0=p(x_{1})p(x_{2})[a(x_{1}),b(x_{2})]
={\bf k}p(x_{1})p(x_{2})\frac{\partial}{\partial x_{2}}
x_{2}^{-1}\delta\left(\frac{x_{1}}{x_{2}}\right).
\end{eqnarray}
Noticing that
\begin{eqnarray}
\Res_{x_{2}}p(x_{1})p(x_{2})\frac{\partial}{\partial x_{2}}
x_{2}^{-1}\delta\left(\frac{x_{1}}{x_{2}}\right)
=-\Res_{x_{2}}p(x_{1})p'(x_{2})x_{2}^{-1}
\delta\left(\frac{x_{1}}{x_{2}}\right)
=-p(x_{1})p'(x_{1}),
\end{eqnarray}
we get ${\bf k}p(x_{1})p'(x_{1})=0$, which implies that ${\bf k}=0$ on $W$.
\end{proof}

In the following we give a slight generalization of
Theorems \ref{tfinite-dimensional-module} and \ref{tcp}.
First recall from [Li2] (Lemmas 2.10 and 2.11) the following result 
(which might be well known, but we do not know any other reference):

\bl{ldensity}
Let $A_{1}$ and $A_{2}$ be associative algebras (with identity)
and let $U_{1}$ and 
$U_{2}$ be irreducible modules for $A_{1}$ and $A_{2}$, respectively.
If either $\End _{A_{1}}U_{1}=\C$, or $A_{1}$ is of countable dimension,
then $U_{1}\otimes U_{2}$ is an irreducible $A_{1}\otimes A_{2}$-module.
\el

The following result slightly generalizes 
the first assertion of Theorem \ref{tfinite-dimensional-module}
of Chari and Presslely with a slightly different proof:

\bp{plevel-0-general}
Assume that $\g$ is of countable dimension.
Let $U_{1},\dots,U_{r}$ 
be irreducible $\g$-modules and let $z_{1},\dots,z_{r}$
be distinct nonzero complex numbers. Then 
$U_{1}(z_{1})\otimes\cdots \otimes U_{r}(z_{r})$ is an irreducible 
$\hat{\g}$-module.
\ep

\begin{proof} Notice that the universal enveloping algebra $U(\hat{\g})$
is of countable dimension.
It follows from Lemma \ref{ldensity}
(and induction) that $U_{1}(z_{1})\otimes \cdots \otimes
U_{r}(z_{r})$ is an irreducible module for the product Lie algebra
$\hat{\g}\oplus \cdots \oplus \hat{\g}$ ($r$ copies).
Denote by $\pi$ the representation homomorphism map.
For $1\le i\le r$, denote by $\psi_{i}$ 
the  $i$-th embedding of $\hat{\g}$ into 
$\hat{\g}\oplus \cdots \oplus \hat{\g}$ ($r$ copies)
and denote by $\psi$ the diagonal map from $\hat{\g}$
to $\hat{\g}\oplus \cdots \oplus \hat{\g}$ ($r$ copies).
Then $\psi =\psi_{1}+\cdots +\psi_{r}.$ 
We also extend the linear maps $\psi$ and $\psi_{1},\dots,\psi_{r}$
on $\hat{\g}[[x,x^{-1}]]$ canonically.

For $1\le i\le r$, set 
$p_{i}(x)=\prod_{j\ne i}(x-z_{j})/(z_{i}-z_{j})$.
Then
\begin{eqnarray}
p_{i}(x)\delta\left(\frac{z_{j}}{x}\right)
=p_{i}(z_{j})\delta\left(\frac{z_{j}}{x}\right)
=\delta_{i,j}\delta\left(\frac{z_{j}}{x}\right)
\end{eqnarray}
for $i,j=1,\dots,r$. Using (\ref{eevaluation-generating})
we have that on $U_{1}(z_{1})\otimes \cdots \otimes U_{r}(z_{r})$,
\begin{eqnarray}
p_{i}(x)\pi \psi_{j}(a(x))=\delta_{i,j}\pi \psi_{j}(a(x))
\;\;\;\mbox{ for }1\le i, j\le r,\; a\in \g.
\end{eqnarray}
Thus on $U_{1}(z_{1})\otimes \cdots \otimes U_{r}(z_{r})$,
\begin{eqnarray}
p_{i}(x)\pi\psi (a(x))=\pi \psi_{i}(a(x))
\;\;\;\mbox{ for }1\le i\le r,\;a\in \g,
\end{eqnarray}
which implies that
\begin{eqnarray}
\pi \psi_{i}(\hat{\g})\subset \pi\psi (\hat{\g})\;\;\;\mbox{ for }i=1,\dots,r.
\end{eqnarray}
{}From this we have
\begin{eqnarray}
\pi \psi (\hat{\g})
=\pi \psi_{1}(\hat{\g})+ \cdots + \pi \psi_{r}(\hat{\g}).
\end{eqnarray}
It follows that $U_{1}(z_{1})\otimes\cdots \otimes U_{r}(z_{r})$ 
is an irreducible $\hat{\g}$-module.
\end{proof}

We also have the following result 
which generalizes Theorem \ref{tcp} of Chari and
Presslely with a different proof: 

\bp{ptensor-product-module}
Assume that $\g$ is of countable dimension.
Let $W$ be an irreducible restricted $\hat{\g}$-module 
(in the category ${\mathcal{R}}$) and 
let $U$ be an irreducible $\hat{\g}$-module in the category $\E$. 
Then the tensor product module $W\otimes U$ is irreducible.
\ep

\begin{proof} Let $M$ be any nonzero
submodule of the tensor product $\hat{\g}$-module
$W\otimes U$. We must prove that $M=W\otimes U$.
Since $W$ and $U$ are irreducible $\hat{\g}$-modules
and $U(\hat{\g})$ is of countable dimension,
by Lemma \ref{ldensity}
$W\otimes U$ is an irreducible $\hat{\g}\oplus \hat{\g}$-module.
Now, it suffices to prove that $M$ is a $\hat{\g}\oplus \hat{\g}$-submodule
of $W\otimes U$ and furthermore it suffices to prove that
\begin{eqnarray}\label{eanm-M}
(a(n)\otimes 1)M\subset M\;\;\;\mbox{ for }a\in \g,\; n\in \Z.
\end{eqnarray}
(Notice that $(1\otimes a(n))w=a(n)w-(a(n)\otimes 1)w$ for $w\in M$.)

By definition there exists a nonzero polynomial
$p(x)$ such that $p(x)a(x)=0$ on $U$ for all $a\in \g$, so that 
\begin{eqnarray}\label{etensor=first}
p(x)(a(x)\otimes 1+1\otimes a(x))=p(x)a(x)\otimes 1
\;\;\mbox{ on }\;W\otimes U. 
\end{eqnarray}
With $M$ being a $\hat{\g}$-submodule of the tensor product module
and with $W$ being a restricted module we have
\begin{eqnarray}
p(x)(a(x)\otimes 1+1\otimes a(x))M\subset M[[x,x^{-1}]],\;\;\;\;
p(x)(a(x)\otimes 1)M\subset (W\otimes U)((x)).
\end{eqnarray}
{}From this, using (\ref{etensor=first}) we have
\begin{eqnarray}
p(x)(a(x)\otimes 1)M\subset M((x))\;\;\;\mbox{ for }a\in \g.
\end{eqnarray}
Let $f(x)$ be the formal Laurent series of rational function $1/p(x)$ 
at zero, so that $f(x)\in \C((x))$.
Then we have
$$a(x)\otimes 1=(f(x)p(x))(a(x)\otimes 1)=f(x)(p(x)(a(x)\otimes 1))$$
on $M$. Consequently,
\begin{eqnarray}
(a(x)\otimes 1)M\subset M((x))\;\;\;\mbox{ for }a\in \g.
\end{eqnarray}
This proves (\ref{eanm-M}), completing the proof.
\end{proof}

The following result tells us when two $\hat{\g}$-modules
of the form $W\otimes U$ obtained in 
Proposition \ref{ptensor-product-module} are isomorphic:

\bp{pidentification-2}
Let $W_{1},W_{2}$ be irreducible $\hat{\g}$-modules in category $\CR$
and let $U_{1}$ and $U_{2}$ be irreducible $\hat{\g}$-modules 
in category $\E$.
Then the tensor product $\hat{\g}$-modules $W_{1}\otimes U_{1}$ 
and $W_{2}\otimes U_{2}$ are isomorphic  if and only if
$W_{1}$ and $U_{1}$ are isomorphic to $W_{2}$ and $U_{2}$, respectively.
\ep

\begin{proof} We only need to prove the ``only if'' part.
Let $f$ be a $\hat{\g}$-module isomorphism from $W_{1}\otimes U_{1}$
onto $W_{2}\otimes U_{2}$. We have
\begin{eqnarray}\label{e2.36}
f(a(x)\otimes 1+1\otimes a(x))v=(a(x)\otimes 1+1\otimes a(x))f(v)
\;\;\;\mbox{ for }a\in \g,\; v\in W_{1}\otimes U_{1}.
\end{eqnarray}
Let $p(x)$ be a nonzero polynomial such that
$$p(x)a(x)U_{1}=0\;\;\mbox{ and }\;\; p(x)a(x)U_{2}=0\;\;\;\mbox{ for }a\in \g.$$
Using this and (\ref{e2.36}) we get
\begin{eqnarray}
p(x)f((a(x)\otimes 1)v)=p(x)(a(x)\otimes 1)f(v)
\;\;\;\mbox{ for }a\in \g,\; v\in W_{1}\otimes U_{1}.
\end{eqnarray}
In view of Remark \ref{rassociativity}, we have
\begin{eqnarray}\label{efav}
f((a(x)\otimes 1)v)=(a(x)\otimes 1)f(v)
\;\;\;\mbox{ for }a\in \g,\; v\in W_{1}\otimes U_{1}.
\end{eqnarray}
Let $u_{2}^{i}$ for $i\in S$ be a basis of $U_{2}$. Then 
$W_{2}\otimes U_{2}=\coprod_{i\in S}W_{2}\otimes \C u_{2}^{i}$.
Denote by $\phi_{i}$ the projection of $W_{2}\otimes U_{2}$ 
onto $W_{2}\otimes \C u_{2}^{i}$.
We have
$$\phi_{i}(a(x)\otimes 1)w=(a(x)\otimes 1)\phi_{i}(w)
\;\;\;\mbox{ for }a\in \g,\; w\in
W_{2}\otimes U_{2},$$
so that
$$\phi_{i}f((a(x)\otimes 1)v)
=(a(x)\otimes 1)\phi_{i}f(v)\;\;\;\mbox{ for }a\in \g,\; v\in
W_{1}\otimes U_{1}.$$
Let $0\ne u_{1}\in U_{1}$. There exists an $i\in S$ 
such that $\phi_{i}f\ne 0$ on 
$W_{1}\otimes \C u_{1}$. We see that the map $\phi_{i}f$ gives rise to a
nonzero $\hat{\g}$-module homomorphism from $W_{1}\;(=W_{1}\otimes \C u_{1})$ 
onto $W_{2}\;(=W_{2}\otimes \C u_{2}^{i})$. 
Because $W_{1}$ and $W_{2}$ are irreducible, 
this nonzero homomorphism is an isomorphism.
This proves that $W_{1}$ is isomorphic to $W_{2}$.

{}From (\ref{e2.36}) and (\ref{efav}) we have
\begin{eqnarray}
f((1\otimes a(x))v)=(1\otimes a(x))f(v)
\;\;\;\mbox{ for }a\in \g,\; v\in W_{1}\otimes U_{1}.
\end{eqnarray}
Then using the same strategy,
we see that $U_{1}$ is isomorphic to $U_{2}$.
\end{proof}

Furthermore, the following result, which is a version of a result of Chari in [C],
gives the equivalence on evaluation $\hat{\g}$-modules
(in category $\E$):

\bp{pisomorphism-E}
Let $U_{1},\dots,U_{r},\; V_{1},\dots,V_{s}$ be nontrivial irreducible
$\g$-modules and let $z_{1},\dots,z_{r}$ and
$\xi_{1},\dots,\xi_{s}$ be two groups of distinct nonzero complex numbers.
Then the $\hat{\g}$-module $U_{1}(z_{1})\otimes \cdots \otimes U_{r}(z_{r})$ is isomorphic to
$V_{1}(\xi_{1})\otimes \cdots \otimes V_{s}(\xi_{s})$ if and only if
$r=s$, $z_{i}=\xi_{i}$ and $U_{i}\cong V_{i}$ up to a permutation.
\ep

\begin{proof} We only need to prove the ``only if'' part.
Let $U$ be any $\hat{\g}$-module in category $\E$.
There exists a (unique nonzero) monic polynomial $p(x)$ 
of least degree such that $p(x)a(x)U=0$ for $a\in \g$. 
Clearly, isomorphic
$\hat{\g}$-modules in category $\E$ have the same monic polynomial.
If $U=U_{1}(z_{1})\otimes \cdots \otimes U_{r}(z_{r})$, we are going to show that
$p(x)=(x-z_{1})\cdots (x-z_{r})$ is the associated monic polynomial.
First, by Lemma \ref{lchari1} we have that $p(x)a(x)=0$
on $U_{1}(z_{1})\otimes \cdots \otimes U_{r}(z_{r})$ for $a\in \g$.
Let $q(x)$ be any polynomial such that $q(x)a(x)=0$ 
on $U_{1}(z_{1})\otimes \cdots \otimes U_{r}(z_{r})$ for $a\in \g$.
Set $p_{i}(x)=\prod_{j\ne i}(x-z_{j})/(z_{i}-z_{j})$ for $i=1,\dots,r$ as
in the proof of Proposition \ref{plevel-0-general}.
For $a\in \g,\; u_{i}\in U_{i}$ with $i=1,\dots,r$, we have
\begin{eqnarray*}
0&=&q(x)p_{i}(x)a(x)(u_{1}\otimes \cdots\otimes u_{r})\nonumber\\
&=&q(x)x^{-1}\delta\left(\frac{z_{i}}{x}\right)
(u_{1}\otimes \cdots \otimes au_{i}\otimes\cdots\otimes u_{r})\nonumber\\
&=&q(z_{i})x^{-1}\delta\left(\frac{z_{i}}{x}\right)
(u_{1}\otimes \cdots \otimes au_{i}\otimes\cdots\otimes u_{r}).
\end{eqnarray*}
Since each $U_{i}$ is a nontrivial $\g$-module,
we must have $q(z_{i})=0$ for $i=1,\dots,r$. Thus $p(x)$ divides $q(x)$.
This proves that $p(x)$ is the associated monic polynomial.

Assume that $U_{1}(z_{1})\otimes \cdots\otimes U_{r}(z_{r})$ is isomorphic
to $V_{1}(\xi_{1})\otimes\cdots\otimes V_{s}(\xi_{s})$ with $f$ a $\hat{\g}$-module
isomorphism map.
Then the two tensor product modules must have the same associated monic polynomial. 
That is, $(x-z_{1})\cdots (x-z_{r})=(x-\xi_{1})\cdots (x-\xi_{s})$.
Thus $r=s$ and up to a permutation $z_{i}=\xi_{i}$ for $i=1,\dots,r$.
Assume that $z_{i}=\xi_{i}$ for $i=1,\dots,r$.
For $1\le i\le r,\; a\in \g$ and for $u_{j}\in U_{j},\; v_{j}\in V_{j}$ 
with $j=1,\dots,r$, we have
\begin{eqnarray}
& &p_{i}(x)a(x)(u_{1}\otimes\cdots\otimes u_{r})
=x^{-1}\delta\left(\frac{z_{i}}{x}\right)
(u_{1}\otimes \cdots \otimes au_{i}\otimes\cdots\otimes u_{r}),\\
& &p_{i}(x)a(x)(v_{1}\otimes\cdots\otimes v_{r})
=x^{-1}\delta\left(\frac{z_{i}}{x}\right)
(v_{1}\otimes \cdots \otimes av_{i}\otimes\cdots\otimes v_{r}).
\end{eqnarray}
Then
\begin{eqnarray}
f(u_{1}\otimes \cdots \otimes au_{i}\otimes\cdots\otimes u_{r})
&=&\Res_{x}x^{-1}\delta\left(\frac{z_{i}}{x}\right)
f(u_{1}\otimes \cdots \otimes au_{i}\otimes\cdots\otimes u_{r})\nonumber\\
&=&\Res_{x}f\left(p_{i}(x)a(x)(u_{1}\otimes\cdots\otimes u_{r})\right)\nonumber\\
&=&\Res_{x}p_{i}(x)a(x)f(u_{1}\otimes\cdots\otimes u_{r})\nonumber\\
&=&\sigma_{i}(a)f(u_{1}\otimes\cdots\otimes u_{r}),
\end{eqnarray}
where for $a\in \g,\; v_{1}\in V_{1},\dots,v_{r}\in V_{r}$,
$$ \sigma_{i}(a)(v_{1}\otimes\cdots\otimes v_{r})=
(v_{1}\otimes \cdots \otimes av_{i}\otimes\cdots\otimes v_{r}). $$
Now, from the proof of Proposition \ref{pidentification-2} we see that
$U_{i}$ is isomorphic to $V_{i}$.
\end{proof}

In view of Remark \ref{runification} and
Proposition \ref{ptensor-product-module}, naturally one wants to know whether
irreducible $\hat{\g}$-modules of the form $W\otimes U$
as in Proposition \ref{ptensor-product-module} exhaust 
the irreducible $\hat{\g}$-modules in the category ${\mathcal{C}}$
up to isomorphism. In the next section
we shall prove that this is true if we restrict ourselves to
integrable module for
a standard affine Lie algebra $\hat{\g}$.

\section{Classification of
irreducible integrable $\hat{\g}$-modules
in the categories $\CR$, $\E$ and ${\mathcal{C}}$}

In this section we  classify irreducible integrable
$\hat{\g}$-modules in the categories $\E$ and ${\mathcal{C}}$
for a standard affine Lie algebra $\hat{\g}$
(with $\g$ a finite-dimensional simple Lie algebra and with $\<\cdot,\cdot\>$ 
the normalized Killing form). It has been proved in [DLM] that
every irreducible integrable $\hat{\g}$-module in the category $\CR$ 
is a highest weight module and every integrable $\hat{\g}$-module 
in the category $\CR$ is completely reducible.
We here show that every irreducible integrable $\hat{\g}$-module
in the category $\E$ is isomorphic to a finite-dimensional evaluation module
and that every irreducible integrable $\hat{\g}$-module
in the category ${\mathcal{C}}$ is isomorphic to
a tensor product of a highest weight integrable module with a
finite-dimensional evaluation module, constructed by Chari and Presslely.

We start with some formal calculus.
First we have ([Li1], [LL])
\begin{eqnarray}\label{edelta=0}
(x_{1}-x_{2})^{m}\left(\frac{\partial}{\partial x_{2}}\right)^{n}
x_{2}^{-1}\delta\left(\frac{x_{1}}{x_{2}}\right)
=0
\end{eqnarray}
for $m>n\ge 0$, and we have
\begin{eqnarray}\label{edelta=delta}
(x_{1}-x_{2})^{m}\frac{1}{n!}\left(\frac{\partial}{\partial x_{2}}\right)^{n}
x_{2}^{-1}\delta\left(\frac{x_{1}}{x_{2}}\right)
=\frac{1}{(n-m)!}\left(\frac{\partial}{\partial x_{2}}\right)^{n-m}
x_{2}^{-1}\delta\left(\frac{x_{1}}{x_{2}}\right)
\end{eqnarray}
for $0\le m\le n$.

\bd{dbarE}
{\em Let $W$ be any vector space. Following [LL] (cf. [Li1]) we set
\begin{eqnarray}
\E(W)=\Hom (W,W((x))).
\end{eqnarray}
We define $\bar{\E}(W)$ to be the subspace of $(\End W)[[x,x^{-1}]]$,
consisting of formal series $a(x)$ such that
$p(x)a(x)\in \Hom (W,W((x)))$ for some nonzero polynomial $p(x)$.
Define $\bar{\E}_{0}(W)$ to be the subspace of $\bar{\E}(W)$
consisting of the formal series $a(x)$ such that
$p(x)a(x)=0$ for some nonzero polynomial $p(x)$.}
\ed

\br{rmodification}
{\em If $a(x)\in \bar{\E}(W)$ and if $x^{m}f(x)a(x)\in \Hom (W,W((x)))$
for some integer $m$ and for some polynomial $f(x)$, then
$f(x)a(x)\in \Hom (W,W((x)))$. In view of this, 
if we need, we may assume that $p(0)\ne 0$ 
for our nonzero polynomial $p(x)$. }
\er

Let $\C(x)$ denote the algebra of rational functions of $x$. 
We define $\iota_{x;0}$ 
to be the linear map from $\C(x)$ to $\C((x))$ such that
for $f(x)\in \C(x)$, $\iota_{x;0}(f(x))$ is the formal Laurent series 
of $f(x)$ at $0$. 
Notice that both $\C(x)$ and $\C((x))$ are
(commutative) fields.
The linear map $\iota_{x;0}$ is a field embedding.
If $p(x)$ is a polynomial
with $p(0)\ne 0$, then $\iota_{x;0}(p(x))\in \C[[x]]$.

\bd{dtildemap}
{\em Let $W$ be a vector space. Define a linear map
$$\psi_{{\mathcal{R}}}: \bar{\E}(W)\rightarrow \E(W)\;(=\Hom (W,W((x))))$$
by
\begin{eqnarray}\label{etildemap}
\psi_{{\mathcal{R}}}(a(x))w=\iota_{x;0}(f(x)^{-1})(f(x)a(x)w)
\;\;\;\mbox{ for }a(x)\in \bar{\E}(W), \;w\in W,
\end{eqnarray}
where $f(x)$ is any nonzero polynomial such that 
$f(x)a(x)\in \Hom (W,W((x)))$.}
\ed

First of all, the map $\psi_{{\mathcal{R}}}$ is well defined; 
the expression on the right-hand side 
of (\ref{etildemap}) makes sense (which is clear)
and does not depend on the choice of $f(x)$. 
Indeed,  let $0\ne f,g\in \C[x]$ be such that 
$$f(x)a(x),\;\; g(x)a(x)\in \Hom (W,W((x))).$$
Set $h(x)=f(x)g(x)$. 
Then $h(x)a(x)\in \Hom(W,W((x)))$. For $w\in W$, we have
\begin{eqnarray*}
\iota_{x;0}(h(x)^{-1})(h(x)a(x)w)
&=&\iota_{x;0}(h(x)^{-1})f(x)(g(x)a(x)w)\\
&=&\iota_{x;0}(g(x)^{-1})(g(x)a(x)w).
\end{eqnarray*}
Similarly, we have
\begin{eqnarray*}
\iota_{x;0}(h(x)^{-1})(h(x)a(x)w)=\iota_{x;0}(f(x)^{-1})(f(x)a(x)w).
\end{eqnarray*}

\br{riotacaancellation}
{\em Note that the expression $\iota_{x;0}(f(x)^{-1})a(x)w$
may not exist in $W[[x,x^{-1}]]$. Thus, in (\ref{etildemap}), 
it is necessary to use all the parenthesis.}
\er

The following is an immediate consequence of (\ref{etildemap})
and the associativity law (recall Remark \ref{rassociativity}):
\bl{lbasic-property}
For $a(x)\in \bar{\E}(W)$, we have
\begin{eqnarray}
f(x)\psi_{\mathcal{R}}(a(x))=f(x)a(x),
\end{eqnarray}
where $f(x)$ is any nonzero polynomial such that 
$f(x)a(x)\in \Hom (W,W((x)))$.
\el

Furthermore we have the following result:

\bp{pdecomposition}
Let $W$ be any vector space. We have
\begin{eqnarray}\label{edecomposition}
\bar{\E}(W)=\E(W)\oplus \bar{\E}_{0}(W).
\end{eqnarray}
Furthermore, the linear map $\psi_{\mathcal{R}}$ from 
$\bar{\E}(W)$ to $\E(W)$, 
defined in Definition \ref{dtildemap}, 
is the projection map of $\bar{\E}(W)$ onto $\E(W)$,
i.e.,
\begin{eqnarray}\label{eprojections}
\psi_{\mathcal{R}}|_{\E(W)}=1\;\;\mbox{ and }
\;\; \psi_{\mathcal{R}}|_{\bar{\E}_{0}(W)}=0.
\end{eqnarray}
\ep

\begin{proof}
Let $a(x)\in \E(W)=\Hom (W,W((x)))$. 
In Definition \ref{dtildemap} we can take $f(x)=1$, 
so that $\psi_{\mathcal{R}}(a(x))w=a(x)w$ for $w\in W$. Thus 
$\psi_{\mathcal{R}}(a(x))=a(x)$.

Now, let $a(x)\in \bar{\E}_{0}(W)$. 
By definition there is a nonzero polynomial $p(x)$
such that $p(x)a(x)=0$ on $W$, 
so that $p(x)a(x)\in \Hom (W,W((x)))$. From definition we have 
$$\psi_{\mathcal{R}}(a(x))w=\iota_{x;0}(p(x)^{-1})(p(x)a(x)w)=0
\;\;\;\mbox{ for }w\in W.$$
Thus $\psi_{\mathcal{R}}(a(x))=0$. 
This proves the property (\ref{eprojections})
and it follows immediately that the sum $\E(W)+ \bar{\E}_{0}(W)$ is
a direct sum.

Let $a(x)\in \bar{\E}(W)$ and let $0\ne f(x)\in \C[x]$ be such that
$f(x)a(x)\in \Hom (W,W((x)))$. In view of Lemma \ref{lbasic-property} we have
$f(x)\psi_{\mathcal{R}}(a(x))=f(x)a(x)$. Then
$f(x)(a(x)-\psi_{\mathcal{R}}(a(x)))=0$, 
which implies that $a(x)-\psi_{\mathcal{R}}(a(x))\in 
\bar{\E}_{0}(W)$. Thus, $a(x)\in \E(W)\oplus \bar{\E}_{0}(W)$. This proves
that $\bar{\E}(W)\subset \E(W)\oplus \bar{\E}_{0}(W)$, from which we have
(\ref{edecomposition}).
\end{proof}

\bd{dprojection-2}
{\em Let $W$ be a vector space. Denote by $\psi_{\E}$ the projection map
of $\bar{\E}(W)$ onto $\E_{0}(W)$ 
with respect to the decomposition (\ref{edecomposition}).
For $a(x)\in \bar{\E}(W)$ we set
\begin{eqnarray}
\tilde{a}(x)&=&\psi_{{\mathcal{R}}}(a(x)),\\
\check{a}(x)&=&\psi_{\E}(a(x))=a(x)-\psi_{{\mathcal{R}}}(a(x))
=a(x)-\tilde{a}(x).
\end{eqnarray}}
\ed

{}From Lemma \ref{lbasic-property} we have
\begin{eqnarray}
f(x)\tilde{a}(x)&=&f(x)a(x),\\
f(x)\check{a}(x)&=&0
\end{eqnarray}
for any nonzero $f(x)\in \C[x]$ such that $f(x)a(x)\in \Hom(W,W((x)))$.

The following result relates the actions of $\psi_{\mathcal{R}}(a(x))$ and
$a(x)$ on $W$:

\bl{lconnection}
For $a(x)\in \bar{\E}(W),\; n\in \Z,\;w\in W$, we have
\begin{eqnarray}
\psi_{{\mathcal{R}}}(a(x))(n)w=\sum_{i=0}^{r}\beta_{i}a(n+i)w
\end{eqnarray}
for some $r\in \N,\;\beta_{1},\dots,\beta_{r}\in \C$,
depending on $a(x),w$ and $n$.
\el

\begin{proof} 
Let $p(x)$ be a polynomial with $p(0)\ne 0$ such that
$p(x)a(x)\in \Hom (W,W((x)))$. 
Then $x^{k}p(x)a(x)w\in W[[x]]$
for some nonnegative integer $k$. Assume that
\begin{eqnarray}
\iota_{x;0}(1/p(x))=\sum_{i\ge 0}\alpha_{i}x^{i}\in \C[[x]].
\end{eqnarray}
Noticing that $\Res_{x}x^{k+m}p(x)a(x)w=0$ for $m\ge 0$,
we have
\begin{eqnarray}
\psi_{{\mathcal{R}}}(a(x))(n)w&=&\Res_{x}x^{n}\psi_{{\mathcal{R}}}(a(x))w
=\Res_{x}x^{n}\iota_{x;0}(1/p(x))(p(x)a(x)w)\nonumber\\
&=&\Res_{x}\sum_{0\le i\le k-n-1}\alpha_{i}x^{n+i}(p(x)a(x)w)\nonumber\\
&=&\Res_{x}\left(\sum_{0\le i\le k-n-1}\alpha_{i}x^{n+i}p(x)a(x)\right)w.
\end{eqnarray}
Then it follows immediately.
\end{proof}

We also have:

\bl{lcommutatorrelation}
Let 
$$a(x), b(x)\in \bar{\E}(W),\;
c_{0}(x), \dots, c_{r}(x)\in (\End W)[[x,x^{-1}]]$$
be such that on $W$,
\begin{eqnarray}\label{ecommutatorold}
[a(x_{1}),b(x_{2})]=\sum_{i=0}^{r}\frac{1}{i!}c_{i}(x_{2})
\left(\frac{\partial}{\partial x_{2}}\right)^{i}
x_{1}^{-1}\delta \left(\frac{x_{2}}{x_{1}}\right).
\end{eqnarray}
Then $c_{0}(x), \dots, c_{r}(x)\in \bar{\E}(W)$ and
\begin{eqnarray}\label{ecommutatornew}
[\tilde{a}(x_{1}),\tilde{b}(x_{2})]
=\sum_{i=0}^{r}\frac{1}{i!}\tilde{c_{i}}(x_{2})
\left(\frac{\partial}{\partial x_{2}}\right)^{i}
x_{1}^{-1}\delta \left(\frac{x_{2}}{x_{1}}\right).
\end{eqnarray}
\el

\begin{proof} 
Using (\ref{edelta=0}) and (\ref{edelta=delta}), and noticing that
$$\Res_{x_{1}}c_{j}(x_{2})\left(\frac{\partial}{\partial x_{2}}\right)^{r}
x_{1}^{-1}\delta \left(\frac{x_{2}}{x_{1}}\right)
=(-1)^{r}\Res_{x_{1}}c_{j}(x_{2})\left(\frac{\partial}{\partial x_{1}}\right)^{r}
x_{1}^{-1}\delta \left(\frac{x_{2}}{x_{1}}\right)=0$$
for $r\ge 1$, we get
\begin{eqnarray}
c_{i}(x_{2})=\Res_{x_{1}}(x_{1}-x_{2})^{i}[a(x_{1}),b(x_{2})].
\end{eqnarray}
Then it is clear that $c_{i}(x)\in \bar{\E}(W)$ for $i=0,\dots,r$, since
$b(x)\in \bar{\E}(W)$.

Let $0\ne f(x)\in \C[x]$ be such that
$$f(x)a(x)=f(x)\tilde{a}(x),\; f(x)b(x)=f(x)\tilde{b}(x),\;
f(x)c_{i}(x)=f(x)\tilde{c_{i}}(x)$$
for $i=0,\dots, r$. Then by 
multiplying both sides of (\ref{ecommutatorold}) by $f(x_{1})f(x_{2})$
we obtain
\begin{eqnarray}
f(x_{1})f(x_{2})[\tilde{a}(x_{1}),\tilde{b}(x_{2})]
=\sum_{i=0}^{r}\frac{1}{i!}f(x_{1})f(x_{2})\tilde{c_{i}}(x_{2})
\left(\frac{\partial}{\partial x_{2}}\right)^{i}
x_{1}^{-1}\delta \left(\frac{x_{2}}{x_{1}}\right).
\end{eqnarray}
Then we may multiply both sides by 
$\iota_{x_{1};0}(f(x_{1})^{-1})\iota_{x_{2};0}(g(x_{2})^{-1})$
to get (\ref{ecommutatornew}).
\end{proof}

The following is the key factorization result:

\bt{tanalyticmodule}
Let $\pi$ be a representation of $\hat{\g}$ on module $W$
in the category ${\mathcal{C}}$. Define linear maps
$\pi_{\mathcal{R}}$ and $\pi_{\mathcal{E}}$ from $\hat{\g}$ 
to $\End W$ in terms of generating functions by
\begin{eqnarray}
& &\pi_{\mathcal{R}}(a(x)+\alpha {\bf k})
=\psi_{\mathcal{R}}(\pi (a(x)))+\alpha \pi ({\bf k}),\\
& &\pi_{\mathcal{E}}(a(x)+\beta {\bf k})=\psi_{\mathcal{E}}(\pi (a(x)))
\end{eqnarray}
for $a\in \g,\; \alpha,\beta\in \C$, where we extend $\pi$ to
$\hat{\g}[[x,x^{-1}]]$ canonically.
Then 
\begin{eqnarray}\label{epi-projections}
\pi =\pi_{\mathcal{R}}+\pi_{\mathcal{E}}
\end{eqnarray}
and the linear map
\begin{eqnarray}
\hat{\g}\oplus \hat{\g}
&\rightarrow& \End W\nonumber\\
(u,v)&\mapsto&
\pi_{\mathcal{R}}(u)+\pi_{\E}(v)
\end{eqnarray}
defines a representation of $\hat{\g}\oplus \hat{\g}$ on $W$.
If $(W,\pi)$ is irreducible,
$W$ is an irreducible $\hat{\g}\oplus \hat{\g}$-module.
Furthermore, $(W,\pi_{\mathcal{R}})$ is a restricted $\hat{\g}$-module 
(in the category ${\mathcal{R}}$) and $(W,\pi_{\E})$ is a
$\hat{\g}$-module in the category $\E$.
\et

\begin{proof} The relation (\ref{epi-projections}) follows from
Proposition \ref{pdecomposition}.
It follows immediately from the 
defining commutator relations (\ref{edef-commutator})
and Lemma \ref{lcommutatorrelation} that
$(W,\pi_{\mathcal{R}})$ is a $\hat{\g}$-module and
it is clear that it is restricted. 
(We view ${\bf k}$ as an element of $\bar{\E}(W)$.)
Consequently, $(W,\pi_{\E})$ is a $\hat{\g}$-module, since
$\pi_{\E}=\pi -\pi_{\CR}$.

Let $0\ne p(x)\in \C[x]$ be such that 
$p(x)\pi(a(x))\in \Hom (W,W((x)))$ for all $a\in \g$.
Then 
$$p(x)\pi_{\CR} (a(x))=p(x)\pi (a(x)),\;\;\;\;p(x)\pi_{\E} (a(x))=0,$$
so that
\begin{eqnarray}
& &p(x)\pi_{\CR} (a(x))=p(x)\psi_{\CR}\pi (a(x))=p(x)\pi (a(x))\\
 & &p(x)\pi_{\E}(a(x))=0
\end{eqnarray}
for $a\in \g$. From this we have that
$(W,\pi_{\E})$ belongs to the category $\E$.

For $a,b\in \g$, using the commutator relations
(\ref{edef-commutator}) and the basic delta-function property we have
\begin{eqnarray}
& &p(x_{1})[\pi_{\CR}(a(x_{1})),\pi_{\E}(b(x_{2}))]\nonumber\\
&=&p(x_{1})[\pi_{\CR}(a(x_{1})),\pi(b(x_{2}))]
-p(x_{1})[\pi_{\CR}(a(x_{1})),\pi_{\CR}(b(x_{2}))]\nonumber\\
&=&p(x_{1})[\pi (a(x_{1})),\pi (b(x_{2}))]
-p(x_{1})[\pi_{\CR}(a(x_{1})),\pi_{\CR}(b(x_{2}))]\nonumber\\
&=&p(x_{1})\pi ([a,b](x_{2}))x_{1}^{-1}\delta\left(\frac{x_{2}}{x_{1}}\right)
+\<a,b\>\pi ({\bf k})p(x_{1})\frac{\partial}{\partial x_{2}}
x_{1}^{-1}\delta\left(\frac{x_{2}}{x_{1}}\right)\nonumber\\
& &-p(x_{1})\pi_{\CR}([a,b](x_{2}))
x_{1}^{-1}\delta\left(\frac{x_{2}}{x_{1}}\right)
-\<a,b\>\pi_{\CR} ({\bf k})p(x_{1})\frac{\partial}{\partial x_{2}}
x_{1}^{-1}\delta\left(\frac{x_{2}}{x_{1}}\right)
\nonumber\\
&=&p(x_{2})\pi ([a,b](x_{2}))x_{1}^{-1}\delta\left(\frac{x_{2}}{x_{1}}\right)
+\<a,b\>\pi ({\bf k})p(x_{1})\frac{\partial}{\partial x_{2}}
x_{1}^{-1}\delta\left(\frac{x_{2}}{x_{1}}\right)\nonumber\\
& &-p(x_{2})\pi_{\CR}([a,b](x_{2}))
x_{1}^{-1}\delta\left(\frac{x_{2}}{x_{1}}\right)
-\<a,b\>\pi ({\bf k})p(x_{1})\frac{\partial}{\partial x_{2}}
x_{1}^{-1}\delta\left(\frac{x_{2}}{x_{1}}\right)
\nonumber\\
&=&0.
\end{eqnarray}
Since $\pi_{\CR}(a(x_{1}))\in \Hom (W,W((x_{1})))$, we can
multiply both sides by $\iota_{x_{1},0}1/p(x_{1})$ and use associativity
to get
\begin{eqnarray}
[\pi_{\CR}(a(x_{1})),\pi_{\E}(b(x_{2}))]=0.
\end{eqnarray}
It follows that $(u,v)\mapsto \pi_{\CR}(u)+\pi_{\E}(v)$
defines a representation of $\hat{\g}\oplus \hat{\g}$ on $W$.
With $\pi=\pi_{\CR}+\pi_{\E}$, it is clear that
if $(W,\pi)$ is irreducible,
$W$ is an irreducible $\hat{\g}\oplus \hat{\g}$-module.
\end{proof}

Furthermore, we  have:

\bp{pidentification-1}
Let $(W_{1},\pi_{1})$ and $(W_{2},\pi_{2})$ be 
$\hat{\g}$-modules in the category ${\mathcal{C}}$ and let $f$ be
a $\hat{\g}$-module homomorphism (isomorphism) from $(W_{1},\pi_{1})$ 
to $(W_{2},\pi_{2})$. Then
$f$ is a $\hat{\g}$-module homomorphism (isomorphism) from
$(W_{1},(\pi_{1})_{\CR})$ to $(W_{2},(\pi_{2})_{\CR})$
and a $\hat{\g}$-module homomorphism (isomorphism) from
$(W_{1},(\pi_{1})_{\E})$ to $(W_{2},(\pi_{2})_{\E})$.
\ep

\begin{proof} Let $p(x)$ be a nonzero polynomial such that for every $a\in \g$,
$$p(x)\pi_{1}(a(x))\in \Hom (W_{1},W_{1}((x))),
\;\; p(x)\pi_{2}(a(x))\in \Hom (W_{2},W_{2}((x))).$$
Then we have
$$p(x)\psi_{\CR}(\pi_{1}(a(x)))=p(x)\pi_{1}(a(x)),\;\;
p(x)\psi_{\CR}(\pi_{2}(a(x)))=p(x)\pi_{2}(a(x)),$$
so that
\begin{eqnarray}
& &p(x)(\pi_{1})_{\CR}(a(x))=
p(x)\psi_{\CR}(\pi_{1}(a(x)))=p(x)\pi_{1}(a(x)),\\
& &p(x)(\pi_{2})_{\CR}(a(x))=
p(x)\psi_{\CR}(\pi_{2}(a(x)))=p(x)\pi_{2}(a(x)).
\end{eqnarray}
For $a\in \g,\; w_{1}\in W_{1}$, we have
\begin{eqnarray}
& &p(x)f((\pi_{1})_{\CR}(a(x))w_{1})=p(x)f(\pi_{1}(a(x))w_{1})
=p(x)\pi_{2}(a(x))f(w_{1})\nonumber\\
& &\hspace{1cm}=p(x)(\pi_{2})_{\CR}(a(x))f(w_{1}).
\end{eqnarray}
Since $f((\pi_{1})_{\CR}(a(x))w_{1}),\;
(\pi_{2})_{\CR}(a(x))f(w_{1})\in W_{2}((x))$, in view of 
Remark \ref{rassociativity} we have
\begin{eqnarray}
f((\pi_{1})_{\CR}(a(x))w_{1})=(\pi_{2})_{\CR}(a(x))f(w_{1})\;\;\;\mbox{ for }a\in \g.
\end{eqnarray}
This proves that $f$ is a $\hat{\g}$-module homomorphism from 
$(W_{1},(\pi_{1})_{\CR})$ to $(W_{2},(\pi_{2})_{\CR})$.
(Notice that $a\otimes t^{n}$ for $a\in \g,\; n\in \Z$ generates
$\hat{\g}$.) Because $(\pi_{i})_{\E}=\pi_{i}-(\pi_{i})_{\CR}$ for $i=1,2$,
it follows that $f$ is also a $\hat{\g}$-module homomorphism from 
$(W_{1},(\pi_{1})_{\E})$ to $(W_{2},(\pi_{2})_{\E})$.
\end{proof}

In view of Theorem \ref{tanalyticmodule},
to classify irreducible representations of $\hat{\g}$
in the category ${\mathcal{C}}$ we need to classify
irreducible representations of $\hat{\g}\oplus \hat{\g}$
which are composed of a representation of $\hat{\g}$
in the category ${\mathcal{R}}$ and a representation of $\hat{\g}$
in the category $\E$. Motivated by this, we next present
some elementary results (or facts) about modules for a tensor product
associative algebra $A_{1}\otimes A_{2}$.

\br{rclassical}
{\em We here collect some facts for general (maybe
infinite-dimensional) associative algebras,
which follow from the proofs for the finite-dimensional case.
The first fact is that if $A$ is an associative algebra (with identity),
$U$ a finitely generated $A$-module and $W=\coprod_{i\in I}W_{i}$
a direct sum of $A$-modules, then
$\Hom _{A}(U,W)\cong \coprod_{i\in I}\Hom_{A}(U,W_{i})$.
With this fact, using the usual proof one can prove the second fact:
Let $A_{1}$ and $A_{2}$ be associative algebras (with identity), let
$W$ be an $A_{1}\otimes A_{2}$-module such that
$W$ viewed as an $A_{1}$-module is completely reducible and 
let $\{U_{i}\;|\;i\in I\}$ be a complete set of representatives 
of equivalence classes of irreducible $A_{1}$-submodules of $W$. 
Assume that $\End_{A_{1}}U_{i}=\C$ for $i\in I$.
Then $W\cong \coprod_{i\in I}U_{i}\otimes \Hom_{A_{1}}(U_{i},W)$,
as an $A_{1}\otimes A_{2}$-module. 
A version of Schur lemma (cf. [Di]) 
is that if $A$ is an associative algebra (with identity) of countable
dimension, then $\End_{A}U=\C$ for any irreducible $A$-module $U$.
In view of this, for the second fact,
the condition that $\End_{A_{1}}U_{i}=\C$ can be replaced by
that condition that $A_{1}$ is of countable dimension.}
\er

The following two lemmas are very useful 
in the proof of our main theorems later:

\bl{ltensor-decomposition-2}
Let $A_{1}$ and $A_{2}$ be associative algebras (with identity)
and let $U$ be an irreducible $A_{1}\otimes A_{2}$-module. Suppose that 
$A_{1}$ is of countable dimension and that 
$U$ as an $A_{1}$-module has an irreducible submodule. Then
$U$ is isomorphic to an $A_{1}\otimes A_{2}$-module of the form
$U_{1}\otimes U_{2}$ as in Lemma \ref{ldensity}.
\el

\begin{proof}
Let $U_{1}$ be an irreducible $A_{1}$-submodule of $U$. Since $U$ is 
an irreducible $A_{1}\otimes A_{2}$-module, we have
$U=(A_{1}\otimes A_{2})U_{1}=A_{2}U_{1}$. For any $a\in A_{2}$, $u\mapsto au$
is an $A_{1}$-homomorphism from $U_{1}$ to $U$. Consequently, for
$a\in A_{2}$, either $aU_{1}=0$ or $aU_{1}$ is an irreducible 
$A_{1}$-submodule isomorphic to $U_{1}$. It follows that $U$
as an $A_{1}$-module is a direct sum of irreducible submodules
isomorphic to $U_{1}$. Furthermore, since $A_{1}$ is of countable dimension, 
{}from Remark \ref{rclassical} we have
$W\cong U_{1}\otimes \Hom_{A_{1}}(U_{1},U)$, where 
$\Hom_{A_{1}}(U_{1},U)$ is a natural $A_{2}$-module which is
necessarily irreducible. 
\end{proof}

\bl{ltensor-decomposition}
Let $A_{1}$ and $A_{2}$ be associative algebras (with identity)
and let $W$ be an $A_{1}\otimes A_{2}$-module.
Assume that $A_{1}$ is of countable dimension and assume that
$W$ is a completely reducible $A_{1}$-module and
a completely reducible $A_{2}$-module.
Then $W$ is isomorphic to a direct sum
of irreducible $A_{1}\otimes A_{2}$-modules of the form
$U\otimes V$ with $U$ an
irreducible $A_{1}$-module and $V$ an irreducible $A_{2}$-module.
\el

\begin{proof} Let $\{ U_{1}^{(i)}\;|\; i\in I\}$ be a
complete set of representatives of equivalence classes of
irreducible $A_{1}$-submodules
of $W$. With $A_{1}$-being countable dimensional,
{}from Remark \ref{rclassical} we have
$$W\cong \coprod_{i\in I}U_{1}^{(i)}\otimes \Hom_{A_{1}}(U_{1}^{(i)},W)$$
as an $A_{1}\otimes A_{2}$-module. Since $W$ is
a completely reducible $A_{2}$-module, $\Hom_{A_{1}}(U_{1}^{(i)},W)$
is a completely reducible $A_{2}$-module. Now it follows from
Lemma \ref{ldensity} that $W$ is a completely reducible
$A_{1}\otimes A_{2}$-module.
\end{proof}

We now classify finite-dimensional irreducible
$\hat{\g}$-modules in category $\E$.
For $a\in \g$, we have (cf. [HL]) 
\begin{eqnarray}
a(x)=\sum_{n\in \Z}(a\otimes t^{n})x^{-n-1}
=a\otimes x^{-1}\delta\left(\frac{t}{x}\right).
\end{eqnarray}
For $f(x)\in \C [x],\; m\in \Z,\; a\in \g$, we have
\begin{eqnarray}
x^{m}f(x)a(x)=a\otimes x^{m}f(x)x^{-1}\delta\left(\frac{t}{x}\right)
=a\otimes t^{m}f(t)x^{-1}\delta\left(\frac{t}{x}\right),
\end{eqnarray}
so that 
\begin{eqnarray}
\Res_{x}x^{m}f(x)a(x)=a\otimes t^{m}f(t).
\end{eqnarray}
It follows immediately that for any $\hat{\g}$-module $W$,
$f(x)a(x)W=0$ if and only if
$(a\otimes f(t)\C[t,t^{-1}])W=0$.
For a nonzero polynomial $p(x)$, we define a subcategory
${\mathcal{E}}_{p}$ of ${\mathcal{E}}$, consisting of
$\hat{\g}$-modules $W$ such that
\begin{eqnarray}
p(x)a(x)w=0\;\;\;\mbox{ for }a\in \g,\; w\in W.
\end{eqnarray}
Then a $\hat{\g}$-module in the category $\E_{p(x)}$ exactly amounts to a
module for the Lie algebra $\g\otimes \C[t,t^{-1}]/p(t)\C[t,t^{-1}]$
(recall Lemma \ref{l-level-E}).

\bl{lEfactorization}
Let $p(x)=x^{k}(x-z_{1})\cdots (x-z_{r})$ with $z_{1},\dots,z_{r}$
distinct nonzero complex numbers and with $k\in \N$.
Then any finite-dimensional irreducible $\hat{\g}$-module 
$W$ in the category $\E_{p(x)}$ is isomorphic to a $\hat{\g}$-module
$U_{1}(z_{1})\otimes \cdots \otimes U_{r}(z_{r})$ for
some finite-dimensional irreducible $\g$-modules $U_{1},\dots,U_{r}$.
\el

\begin{proof} Noticing that 
$\C[t,t^{-1}]/p(t)\C[t,t^{-1}]
=\prod_{i=1}^{r}\C[t,t^{-1}]/(t-z_{i})\C[t,t^{-1}]$, we have
\begin{eqnarray}
\g\otimes (\C[t,t^{-1}]/p(t)\C[t,t^{-1}])
=\prod_{i=1}^{r}\g\otimes \C[t,t^{-1}]/(t-z_{i})\C[t,t^{-1}].
\end{eqnarray}
Notice that for any nonzero complex number $z$,
a $\g\otimes \C[t,t^{-1}]/(t-z)\C[t,t^{-1}]$-module exactly amounts to
an evaluation $\hat{\g}$-module $U(z)$.
Set 
$$A_{i}=U\left(\g\otimes \C[t,t^{-1}]/(t-z_{i})\C[t,t^{-1}]\right)$$
for $i=1,\dots,r$.
Since $W$ is finite-dimensional, $W$ viewed as an
$A_{i}$-module contains an irreducible submodule.
It now follows from Lemma \ref{ltensor-decomposition-2} (and induction).
\end{proof}

We also have the following result:

\bp{pE-irreducible-module}
Assume that $[\g,\g]=\g$.
Then any finite-dimensional irreducible $\hat{\g}$-module 
$W$ in the category $\E$ is isomorphic to a $\hat{\g}$-module
$U_{1}(z_{1})\otimes \cdots \otimes U_{r}(z_{r})$ for
some finite-dimensional $\g$-modules $U_{1},\dots,U_{r}$ and
for some distinct nonzero complex numbers $z_{1},\dots,z_{r}$.
\ep

\begin{proof} In view of Lemma \ref{lEfactorization},
it suffices to prove that $W$ is in the category $\E_{p(x)}$
with $p(x)$ a nonzero polynomial whose any nonzero
root is multiplicity-free. In view of Remark \ref{rmodification},
there exists a polynomial 
$p(x)$ with $p(0)\ne 0$ such that $p(x)a(x)W=0$ for $a\in \g$.
Let $p(x)$ be such a monic polynomial with the least degree.
Thus
\begin{eqnarray}
p(x)=(x-z_{1})^{k_{1}}\cdots (x-z_{r})^{k_{r}},
\end{eqnarray}
where $z_{1},\dots,z_{r}$ are distinct nonzero complex numbers 
and $k_{1},\dots,k_{r}$ are positive integers.

Let $I$ be the annihilating ideal of $W$ in $\hat{\g}$.
Then $(\g\otimes p(t)\C[t,t^{-1}])\subset I$ and
$W$ is an irreducible faithful $\hat{\g}/I$-module.
Therefore (cf. [H]) $\hat{\g}/I$ is reductive 
(where we using the fact that $W$ is finite-dimensional).
Set $f(x)=(x-z_{1})\cdots (x-z_{r})$ and
let $k$ be the largest one among $k_{1},\dots,k_{r}$.
We see that $p(x)$ is a factor of $f(x)^{k}$.
It follows that the quotient space $(\g\otimes f(t)\C[t,t^{-1}])/I$ is a solvable ideal of
$\hat{\g}/I$. With $\hat{\g}/I$ being reductive,
$(\g\otimes f(t)\C[t,t^{-1}])/I$ must be in the center of $\hat{\g}/I$.
{}From this we have that $[\g,\g]\otimes f(t)\C[t,t^{-1}]\subset I$,
which implies that
$\g\otimes f(t)\C[t,t^{-1}]\subset I$, since $\g=[\g,\g]$ by assumption.
This proves that $f(x)a(x)W=0$ for $a\in \g$. Consequently,
$f(x)=p(x)$, that is, $k_{1}=\cdots =k_{r}=1$.
\end{proof}

\br{rcomments}
{\em  In Proposition \ref{pE-irreducible-module}, the condition 
$\g=[\g,\g]$ is necessary. For example, let $\g$ be an abelian Lie algebra.
For any nonzero linear functional $\chi$ on $\hat{\g}$ 
with $\psi({\bf k})=0$, we have a
one-dimensional irreducible $\hat{\g}$-module $\C$ 
with $\hat{\g}$ acting according to $\chi$.
In general, such a module may not be in category $\E$. }
\er

For the rest of this section we assume that
{\em $\hat{\g}$ is a standard affine Lie algebra with
$\g$ a finite-dimensional simple Lie algebra and with $\<\cdot,\cdot\>$ 
the normalized Killing form.} We retain all the notations and definitions
in Remark \ref{rintegrable-module}.

The following result is a refinement of
Theorem \ref{tanalyticmodule}:

\bp{pintegrability}
Let $\pi$ be a representation of $\hat{\g}$
on integrable $\hat{\g}$-module $W$ in the category of
${\mathcal{C}}$. Then $(W,\pi_{\mathcal{R}})$ is a
restricted integrable $\hat{\g}$-module (in the category
${\mathcal{R}}$) and $(W,\pi_{\E})$ is 
an integrable $\hat{\g}$-module in the category $\E$.
\ep

\begin{proof} In view of Theorem \ref{tanalyticmodule},
we only need to show that $(W,\pi_{\mathcal{R}})$ 
and $(W,\pi_{\E})$ are integrable $\hat{\g}$-modules.
We must prove that for $a\in \g_{\alpha}$ with $\alpha\in
\Delta$ and for $n\in \Z$, $\tilde{a}(n)$ and $\check{a}(n)$ 
act locally nilpotently on $W$.

Let $a\in \g_{\alpha}$ with $\Delta$ and $n\in \Z$. 
Notice that $[a(r),a(s)]=0$ for $r,s\in \Z$, since $[a,a]=0$ and
$\<a,a\>=0$. For $w\in W$, we have
$$a(r)\tilde{a}(x)w=a(r)\iota_{x;0}(1/p(x)) (p(x)a(x)w)
=\iota_{x;0}(1/p(x)) (p(x)a(x)a(r)w)=\tilde{a}(x)a(r)w.$$
Thus
\begin{eqnarray}\label{ecommuting-tilde}
a(r)\tilde{a}(s)=\tilde{a}(s)a(r)\;\;\;\mbox{ for }r,s\in \Z.
\end{eqnarray}
Let $w\in W$ be an arbitrarily fixed vector.
By Lemma \ref{lconnection}, 
$$\tilde{a}(n)w=\sum_{i=0}^{r}\beta_{i}a(n+i)w$$
for some positive integer $r$ and for some complex numbers
$\beta_{1},\dots,\beta_{r}$.
Using (\ref{ecommuting-tilde}) we get
\begin{eqnarray}\label{eexp-formula}
\tilde{a}(n)^{p}w=(\beta_{0}a(n)+\cdots +\beta_{r}a(n+r))^{p}w
\;\;\;\mbox{ for any }p\ge 0.
\end{eqnarray}
Since $(W,\pi)$ is an integrable $\hat{\g}$-module, there is a
positive integer $k$ such that
$$a(m)^{k}w=0\;\;\;\mbox{ for }m=n,n+1,\dots,n+r.$$
Combining this with (\ref{eexp-formula}) 
we obtain $\tilde{a}(n)^{k(r+1)}w=0$.

Since $\check{a}(n)=a(n)-\tilde{a}(n)$ and
$[a(n),\tilde{a}(n)]=0$, we get
\begin{eqnarray}
\check{a}(n)^{k(r+2)}w=(a(n)-\tilde{a}(n))^{k(r+2)}w
=\sum_{i\ge 0}\binom{k(r+2)}{i}(-1)^{i}a(n)^{k(r+2)-i}
\tilde{a}(n)^{i}w=0.
\end{eqnarray}
This proves that $\tilde{a}(n)$ and $\check{a}(n)$ 
act locally nilpotently on $W$, completing
the proof.
\end{proof}

We shall need the following fact which is a reformulation
of Lemma 3.6 of [DLM]:

\bl{ldlm-basis}
There is a basis $\{ a_{1},\dots,a_{r}\}$ of $\g$ such that
\begin{eqnarray}
[a_{i}(m),a_{i}(n)]=0\;\;\;\mbox{ for }1\le i\le r,\; m,n\in \Z
\end{eqnarray}
and such that
for $1\le i\le r$ and for any $n\in \Z$, $a_{i}(n)$ acts locally nilpotently
on all integrable $\hat{\g}$-modules.
\el

\begin{proof}
For $\alpha\in\Delta_{+}$, choose nonzero vectors
$e_{\alpha}\in \g_{\alpha}, f_{\alpha}\in \g_{-\alpha}, h_{\alpha}\in \h$ 
such that $[h_{\alpha},e_{\alpha}]=2e_{\alpha}$, 
$[h_{\alpha},f_{\alpha}]=-2f_{\alpha}$ and
$[e_{\alpha},f_{\alpha}]=h_{\alpha}$.
Set 
$\sigma_{\alpha}=e^{{\rm ad}e_{\alpha}}$,
an inner automorphism of Lie algebra $\g$.
Then
$\sigma_{\alpha}(f_{\alpha})=f_{\alpha}+h_{\alpha}-e_{\alpha}$. Since 
$\{e_{\alpha}, f_{\alpha}, h_{\alpha}\;|\;\alpha\in \Delta_{+}\}$ 
is a basis of $\g$, 
$\{ e_{\alpha}, f_{\alpha}, \sigma_{\alpha}(f_{\alpha})\;|\;
\alpha\in \Delta_{+}\}$ is also a basis of $\g$. 
On any integrable $\hat{\g}$-module we have (cf. [H], [K1])
\begin{eqnarray}
\exp (e_{\alpha}(0))f_{\alpha}(n)\exp (-e_{\alpha}(0))
=\sigma_{\alpha}(f_{\alpha})(n)\;\;\;\mbox{ for }n\in \Z.
\end{eqnarray}
Since for $n\in \Z$, $f_{\alpha}(n)$ acts locally nilpotently
on any integrable $\hat{\g}$-module, $\sigma_{\alpha}(f_{\alpha})(n)$
also acts locally nilpotently
on any integrable $\hat{\g}$-module.
Then $\{ e_{\alpha},f_{\alpha}, \sigma_{\alpha}(e_{\alpha})\;|\; 
\alpha\in \Delta_{+}\}$,
is a basis of $\g$, satisfying the desired property.
\end{proof}

The following result is a reformulation of Theorem 3.7 
(cf. Remark 3.9) of [DLM]:

\bt{tdlm}
Every nonzero restricted integrable $\hat{\g}$-module is 
a direct sum of (irreducible) highest weight integrable modules. 
In particular, every irreducible integrable $\hat{\g}$-module
$W$ is a highest weight integrable module. 
\et

\begin{proof} 
As in [DLM], in view of the complete reducibility theorem
in [K1] we only need to show that every nonzero restricted integrable 
$\hat{\g}$-module 
$W$ contains a highest weight integrable (irreducible) submodule.
We now reformulate the proof of [DLM, Theorem 3.7] as follows:

{\bf Claim 1:} {\em There exists a nonzero $u\in W$ such that 
$(\g \otimes t{\C}[t])u=0$.} For $n\in \Z$, set
$\g(n)=\{ a(n)\;|\; a\in \g\}$.
For any nonzero $u\in W$, since $W$ is restricted,
$\g(n)u=0$ for $n$ sufficiently large, 
so that $\sum_{n\ge 1}\g(n)u$ is finite-dimensional.
For any $u\in W$, we define $d(u)=\dim \sum_{n\ge 1}\g(n)u$.
If there is a $0\ne u\in W$ such that $d(u)=0$, then
$(\g\otimes t{\C}[t])u=0$. 

Suppose that $d(u)>0$ for any $0\ne u\in W$. 
Take $0\ne u\in W$ such that $d(u)$ is minimal. 
By Lemma \ref{ldlm-basis},
there exists a basis $\{ a_{1},\dots a_{r}\}$ of $\g$ such that
$a_{i}(n)$ locally nilpotently act on $W$ for
$i=1,\dots,r,\; n\in \Z$.
Let $k$ be the positive integer such that
$\g(k)u\ne 0$ and $\g(n)u=0$ whenever $n>k$. By
the definition of $k$, $a_{i}(k)u\ne 0$ for some $1\le i\le r$. 

Notice that $a_{i}(k)^{s}u=0$ for some nonnegative integer $s$.
Let $m$ be the nonnegative integer such that $a_{i}(k)^{m}u\ne 0$ 
and $a_{i}(k)^{m+1}u=0$. 
Set $v=a_{i}(k)^{m}u$. We will obtain
a contradiction by showing that
$d(v)<d(u)$.  
First we prove that if $a(n)u=0$ for some $a\in \g,\; n\ge 1$, 
then $a(n)v=0$. In the following we will show by induction
on $m$ that $a(n)a_{i}(k)^mu=0$ for 
any $a\in \g$ and $m\ge 0$. If $m=0$ this is immediate by the
choice of $u.$ Now assume that the result holds for $m.$  
Since $[a,a_{i}](k+n)u=0$ (from the definition of $k$) and $a(n)u=0,$ 
by the induction assumption that $a(n)a_{i}(k)^mu=0$ we have
\begin{eqnarray}
[a,a_{i}](k+n)a_{i}(k)^{m}u=0,\;a(n)a_{i}(k)^{m}u=0.
\end{eqnarray}
Thus
\begin{eqnarray}
& &a(n)a_{i}(k)^{m+1}u
=[a(n),a_{i}(k)]a_{i}(k)^{m}u+a_{i}(k)a(n)a_{i}(k)^{m}u\nonumber\\
& &\ \ \ \ =[a,a_{i}](k+n)a_{i}(k)^{m}u+a_{i}(k)a(n)a_{i}(k)^{m}u
\nonumber\\
& &\ \ \ \ =0,
\end{eqnarray}
as required. In particular, we see that
$a(n)v=a(n)a_{i}(k)^{r}u=0$.
Therefore, $d(v)\le d(u)$. Since $a_{i}(k)v=0$ and $a_{i}(k)u\ne 0$, we have
$d(v)< d(u)$, a contradiction.

{\bf Claim 2:} $W$ contains an irreducible
highest weight integrable submodule. Set
\begin{eqnarray}
\Omega(W)=\{u\in W\;|\; (\g\otimes t{\C}[t])u=0\}.
\end{eqnarray}
Then $\Omega(W)$ is a $\g$-submodule of $W$ and 
it is nonzero by Claim 1. 
Since $a_{i}(0)$ for $i=1,\dots,r$ act locally nilpotently on
$\Omega(W)$, it follows from the PBW theorem that for any $u\in
\Omega(W)$, $U(\g)u$ 
is finite-dimensional, so that $U(\g)u$ is a direct sum of
finite-dimensional irreducible $\g$-modules.
Let $u\in \Omega(W)$ be a highest weight vector for $\g$.
It is clear that $u$ is a singular vector for $\hat{\g}$.
It follows from [K1] that 
$u$ generates an irreducible $\hat{\g}$-module. 
\end{proof}

We also have the following result 
(cf. Theorem \ref{tfinite-dimensional-module}):

\bp{pevaluation-simple-modules}
The irreducible integrable $\hat{\g}$-modules 
in the category ${\mathcal{E}}$ up to isomorphism
are exactly those evaluation modules $U_{1}(z_{1})\otimes \cdots
\otimes U_{r}(z_{r})$ where $U_{i}$ are finite-dimensional irreducible
$\g$-modules and $z_{i}$ are distinct nonzero complex numbers.
\ep

\begin{proof} In view of Proposition \ref{pE-irreducible-module} 
it suffices to prove that
every irreducible integrable $\hat{\g}$-module $W$ 
in the category ${\mathcal{E}}$ is finite-dimensional.
Since $W$ is in the category $\E$, there is a nonzero polynomial
$p(x)$ such that $(a\otimes p(t)\C[t,t^{-1}])W=0$ for $a\in \g$.
Let $I$ be the annihilating ideal of $W$ in $\hat{\g}$.
Then $\hat{\g}/I$ is finite-dimensional.
Recall from Lemma \ref{ldlm-basis} that there is a basis 
$\{ a_{1},\dots,a_{r}\}$ of $\g$ such that 
for any $1\le i\le r,\; n\in \Z$, $a_{i}(n)$ acts locally nilpotently on $W$.
Let $0\ne w\in W$. Since $W$ is irreducible, we have
$W=U(\hat{\g})w=U(\hat{\g}/I)w$.
In view of the PBW theorem (for $\hat{\g}/I$ using a basis consisting
of the cosets of finitely many $a_{i}(n)$'s) we have
that $W$ is finite-dimensional, completing the proof.
\end{proof}

Now, we are in a position to prove our main result:

\bt{tclassification-simple}
Every irreducible integrable $\hat{\g}$-module 
in the category ${\mathcal{C}}$ 
is isomorphic to a module of the form 
$W\otimes U_{1}(z_{1})\otimes \cdots \otimes U_{r}(z_{r})$,
where $W$ is an irreducible integrable highest weight
$\hat{\g}$-module and $U_{1},\dots,U_{r}$ are finite-dimensional
irreducible $\g$-modules with $z_{1},\dots,z_{r}$ distinct nonzero
complex numbers.
\et

\begin{proof}
Let $\pi$ be an irreducible integrable representation of $\hat{\g}$
on module $W$ in the category ${\mathcal{C}}$.
By Theorem \ref{tanalyticmodule}, $W$ is an irreducible 
$\hat{\g}\oplus \hat{\g}$-module with $(u,v)$ acting as
$\pi_{\mathcal{R}}(u)+\pi_{\mathcal{E}}(v)$ 
for $u,v\in \hat{\g}$ and we have $\pi=\pi_{\mathcal{R}}+\pi_{\mathcal{E}}$.
Furthermore, by Proposition \ref{pintegrability},
$(W,\pi_{\mathcal{R}})$ is an integrable restricted $\hat{\g}$-module
and $(W,\pi_{\E})$
is an integrable $\hat{\g}$-module in the category $\E$. 
In view of Theorem \ref{tdlm}, $(W,\pi_{\mathcal{R}})$ is a direct sum
of integrable highest weight (irreducible) $\hat{\g}$-modules.
Now it follows immediately from Lemma \ref{ltensor-decomposition-2}
with $A_{1}=A_{2}=U(\hat{\g})$ (which is of countable dimension)
and Proposition \ref{pevaluation-simple-modules}.
\end{proof}

Recall (Theorem \ref{tdlm}) that 
every integrable $\hat{\g}$-module in the category $\CR$ 
is completely reducible. But,
an integrable $\hat{\g}$-module in the category 
$\E$ is not necessarily completely reducible.
(Notice that any finite-dimensional $\hat{\g}$-module
in the category $\E$ is integrable, but it is not necessarily 
completely reducible.) Nevertheless we have:

\bp{pcase1}
Let $p(x)$ be a nonzero polynomial such that 
all the nonzero roots are multiplicity-free. 
Then every integrable $\hat{\g}$-module in the category ${\cal{E}}_{p(x)}$ 
is semisimple and every integrable $\hat{\g}$-module 
in the category ${\mathcal{C}}_{p(x)}$ is semisimple.
\ep

\begin{proof} Set $p(x)=x^{k}(x-z_{1})\cdots (x-z_{r})$,
where $k\in \N$ and $z_{1},\dots,z_{r}$ are distinct nonzero complex
numbers. From the proof of Lemma \ref{lEfactorization},
a $\hat{\g}$-module in the category $\E_{p(x)}$ amounts to
a module for the product Lie algebra
$$\coprod_{i=1}^{r}(\g\otimes \C[t,t^{-1}]/(t-z_{i})\C[t,t^{-1}]).$$
Let $W$ be an integrable $\hat{\g}$-module 
in the category ${\cal{E}}_{p(x)}$.
Using the basis of $\g$ as in the proof of
Proposition \ref{pevaluation-simple-modules}
it follows from the PBW theorem that 
any vector in $W$ generates a finite-dimensional
$\hat{\g}$-submodule. 
With $\g\otimes \C[t,t^{-1}]/(t-z_{i})\C[t,t^{-1}]=\g$,
it follows that $W$ as a module for each of the Lie algebras
$\g\otimes \C[t,t^{-1}]/(t-z_{i})\C[t,t^{-1}]$ is completely reducible. 
Now it follows from Lemma \ref{ltensor-decomposition}
that $W$ is completely reducible.

Finally, with the first assertion and Theorem \ref{tdlm}
it follows from Lemma \ref{ltensor-decomposition}
that every integrable $\hat{\g}$-module in the category
${\mathcal{C}}$ is completely reducible.
\end{proof}

\section{A relation between tensor product module $W\otimes U(z)$
and fusion rules}

In this section we relate the tensor product module $W\otimes U(z)$
in the category ${\mathcal{C}}$ with the fusion rules of certain type
for the vertex operator algebra associated with the affine Lie algebra
$\hat{\g}$ of level $\ell$.

As in Section 2, let $\g$ be a (not necessarily finite-dimensional) 
Lie algebra equipped with a nondegenerate symmetric invariant bilinear form
$\<\cdot,\cdot\>$ and let $\hat{\g}$ be the associated affine Lie algebra.
Recall the extended affine Lie algebra (cf. [K1])
\begin{eqnarray}
\tilde{\g}=\hat{\g}\oplus \C {\bf d}
=\g\otimes \C[t,t^{-1}]\oplus \C {\bf k}\oplus \C {\bf d},
\end{eqnarray}
where $[{\bf d},{\bf k}]=0$ and 
\begin{eqnarray}
[{\bf d}, a\otimes t^{n}]
=n(a\otimes t^{n})\;\;\;\mbox{ for }a\in \g,\; n\in \Z.
\end{eqnarray}

A $\tilde{\g}$-module $W$ is said to be {\em upper truncated}
if $W=\coprod_{\lambda\in \C}W(\lambda)$, where for $\lambda\in \C$,
$$W(\lambda)=\{ w\in W\;|\; {\bf d}w=\lambda w\},$$
such that for any $\lambda\in \C$, $W(\lambda+n)=0$ for $\in \Z$ sufficiently large.
Clearly, we have
\begin{eqnarray}
a(n)W(\lambda)\subset W(\lambda+n)\;\;\;\mbox{ for }a\in \g,\; n\in \Z,\; \lambda\in \C.
\end{eqnarray}
Then every upper truncated $\tilde{\g}$-module is a restricted $\hat{\g}$-module.
For an upper truncated $\tilde{\g}$-module $W=\coprod_{\lambda\in \C}W(\lambda)$,
we set (cf. [HL])
\begin{eqnarray}
\overline{W}=\prod_{\lambda\in \C}W(\lambda),
\end{eqnarray}
the formal completion of $W$. Then $\overline{W}$ is again a
$\tilde{\g}$-module (but not a restricted module).

For any $\g$-module $U$, $L(U)=U\otimes \C[t,t^{-1}]$ is naturally a
$\tilde{\g}$-module where
\begin{eqnarray}
& &a(m)(u\otimes t^{n})=au\otimes t^{m+n},\\
& &{\bf d}(u\otimes t^{n})=(n+1)(u\otimes t^{n})
\;\;\;\mbox{ for }a\in \g,\; u\in U,\; m,n\in \Z
\end{eqnarray}
and ${\bf k}$ acts as zero (cf. [CP2], [K1]). 
Such a $\tilde{\g}$-module is often called a loop module.
We have $L(U)=\coprod_{n\in \Z}L(U)(n)$, where
$L(U)(n)=(U\otimes \C t^{n-1})$ for $n\in \Z$.

\br{revalution-modules}
{\em Notice that if $U$ is not a trivial $\g$-module, i.e.,
$\g U=0$, then the action of
$\hat{\g}$ on the evaluation $\hat{\g}$-module $U(z)$ 
cannot be extended to a module action for the extended 
affine Lie algebra for $a\in \g,\; m\in \Z,\; u\in U$, 
$\tilde{\g}$. Otherwise,  we have
\begin{eqnarray}
0=({\bf d}a(m)-a(m){\bf d}-ma(m))u
=z^{m}({\bf d}au-a{\bf d}u-mau)
\end{eqnarray}
for $a\in \g,\; m\in \Z,\; u\in U$, which implies that $au=0$ for
$a\in \g,\; u\in U$, a contradiction. }
\er

Let $W_{1}$ and $W$ be upper truncated
$\tilde{\g}$-modules of level $\ell$ and $U(z)$ be an
evaluation $\hat{\g}$-module (of level zero),
where $z$ is a fixed nonzero complex number.
We have a (tensor product) $\hat{\g}$-module $W_{1}\otimes U(z)$
and a (tensor product) $\hat{\g}$-module $W_{1}\otimes L(U)$.
For homogeneous vector $w_{1}\in W_{1}$ and 
for $u\in U,\; n\in \Z$, we have
\begin{eqnarray}
\deg (w_{1}\otimes u\otimes t^{n})=\deg w_{1}+n+1.
\end{eqnarray}

We next show that there is a canonical linear isomorphism
from $\Hom _{\tilde{\g}}(W_{1}\otimes U\otimes \C[t,t^{-1}],W)$
to $\Hom _{\hat{\g}}(W_{1}\otimes U(z),\overline{W})$.

Let $\psi$ be a $\tilde{\g}$-module homomorphism from
$W_{1}\otimes U\otimes \C[t,t^{-1}]$ to $W$. We define a linear map
\begin{eqnarray}
\hat{\psi}: W_{1}\otimes U&\rightarrow& \overline{W}\nonumber\\
         w_{1}\otimes u&\mapsto& 
\sum_{n\in \Z} z^{-n-1}\psi(w_{1}\otimes u\otimes t^{n}).
\end{eqnarray}
We are going to show that $\hat{\psi}$ is in fact 
a $\hat{\g}$-module homomorphism from
the tensor product module $W_{1}\otimes U(z)$ to $\overline{W}$.

Let $a\in \g,\; m\in \Z,\; w_{1}\in W_{1},\; u\in U$. We have
\begin{eqnarray}
& &\hat{\psi}(a(m) (w_{1}\otimes u))\nonumber\\
&=&\hat{\psi}(a(m)w_{1}\otimes u+w_{1}\otimes z^{m}au)\nonumber\\
&=&\sum_{n\in \Z}z^{-n-1}\psi(a(m)w_{1}\otimes u\otimes t^{n})
+z^{m-n-1}\psi (w_{1}\otimes au\otimes t^{n})\nonumber\\
&=&\sum_{n\in \Z}z^{-n-1}\left(\psi (a(m)w_{1}\otimes u\otimes t^{n})
+\psi(w_{1}\otimes au\otimes t^{m+n})\right)\nonumber\\
&=&\sum_{n\in \Z}z^{-n-1}
\psi \left(a(m)(w_{1}\otimes u\otimes t^{n})\right)\nonumber\\
&=&\sum_{n\in \Z}z^{-n-1}a(m)\psi (w_{1}\otimes u\otimes t^{n})\nonumber\\
&=&a(m)\hat{\psi }(w_{1}\otimes u).
\end{eqnarray}
Since $a\otimes t^{m}$ for $a\in \g,\; m\in \Z$ generate $\hat{\g}$,
$\hat{\psi}$ is a $\hat{\g}$-module homomorphism.
Clearly, $\hat{\psi}=0$ implies $\psi=0$.
Then we obtain a one-to-one linear map 
{}from $\Hom _{\tilde{\g}}(W_{1}\otimes U\otimes \C[t,t^{-1}],W)$
to $\Hom _{\hat{\g}}(W_{1}\otimes U(z),\overline{W})$
sending $\psi$ to $\hat{\psi}$.

On the other hand, let $\phi$ be a $\hat{\g}$-module homomorphism from
$W_{1}\otimes U(z)$ to $\overline{W}$. 
For any $\lambda\in \C$, denote by $p_{\lambda}^{W}$ the projection of
$\overline{W}$ onto the homogeneous subspace $W(\lambda)$.
We have
\begin{eqnarray}
p_{\lambda}^{W}(a(m)\bar{w})=a(m)p_{\lambda-m}^{W}(\bar{w})
\;\;\;\mbox{ for }a\in \g,\; \lambda\in \C,\; m\in \Z,\; 
\bar{w}\in \overline{W}.
\end{eqnarray}
Define a linear map
$\tilde{\phi}$ from $W_{1}\otimes U\otimes \C[t,t^{-1}]$ to $W$ 
by
\begin{eqnarray}
\tilde{\phi}(w_{1}\otimes u\otimes t^{n})
=z^{n+1}p^{W}_{\deg w_{1}+n+1} \phi(w_{1}\otimes u)
\end{eqnarray}
for homogeneous vector $w_{1}\in W_{1}$ and for $u\in U,\; n\in \Z$.
We now show that $\tilde{\phi}$ is a $\tilde{\g}$-module homomorphism.
Let $w_{1}\in W_{1}$ be homogeneous and let $a\in \g,\; u\in U,\; m,n\in \Z$. 
Noticing that $\deg a(m)w_{1}=\deg w_{1}+m$, we have
\begin{eqnarray}
& &\tilde{\phi}(a(m)(w_{1}\otimes u\otimes t^{n}))\nonumber\\
&=&\tilde{\phi}(a(m)w_{1}\otimes u\otimes t^{n}
+ w_{1}\otimes au\otimes t^{m+n})\nonumber\\
&=&z^{n+1}p^{W}_{\deg w_{1}+m+n+1}\phi (a(m)w_{1}\otimes u)
+z^{m+n+1}p^{W}_{\deg w_{1}+m+n+1}\phi (w_{1}\otimes au)\nonumber\\
&=&z^{n+1}p^{W}_{\deg w_{1}+m+n+1}\phi \left(a(m)(w_{1}\otimes u)\right)
\nonumber\\
&=&z^{n+1}p^{W}_{\deg w_{1}+m+n+1}a(m)\phi (w_{1}\otimes u)\nonumber\\
&=&z^{n+1}a(m)p^{W}_{\deg w_{1}+n+1}\phi (w_{1}\otimes u)\nonumber\\
&=&a(m)\tilde{\phi}(w_{1}\otimes u\otimes t^{n}).
\end{eqnarray}
This shows that $\tilde{\phi}$ is indeed a $\tilde{\g}$-module homomorphism.

For $\phi\in \Hom_{\hat{\g}}(W_{1}\otimes U(z),W)$,
set $\psi=\tilde{\phi}$. For homogeneous vector $w_{1}\in W_{1}$ 
and for $u\in U$, we have
\begin{eqnarray}
\hat{\psi}(w_{1}\otimes u)
&=&\sum_{n\in \Z}z^{-n-1}\psi(w_{1}\otimes u\otimes t^{n})\nonumber\\
&=&\sum_{n\in \Z}z^{-n-1}\tilde{\phi}(w_{1}\otimes u\otimes t^{n})\nonumber\\
&=&\sum_{n\in \Z}p^{W}_{\deg w_{1}+n+1}\phi(w_{1}\otimes u)\nonumber\\
&=&\phi(w_{1}\otimes u).
\end{eqnarray}
This shows that the linear map $\psi\mapsto \hat{\psi}$ is 
also onto.

To summarize we have:

\bp{pequivalence}
Let $W_{1},W$ be upper truncated $\tilde{\g}$-modules of level $\ell$ and let $U$
be a $\g$-module and $z$ a nonzero complex number.
Then the map $\psi\mapsto \hat{\psi}$ from 
$\Hom_{\tilde{\g}}(W_{1}\otimes U\otimes \C[t,t^{-1}],W)$ 
to $\Hom_{\hat{\g}}(W_{1}\otimes U(z),\overline{W})$
is a linear isomorphism. The inverse map is
given by $\phi\mapsto \tilde{\phi}$.
\ep

Let $\ell$ be any complex number.
Take $U$ to be the one-dimensional trivial $\g$-module $\C$
in the (\ref{einduced-module}) and set
\begin{eqnarray}
V_{\hat{\g}}(\ell,0)=M_{\hat{\g}}(\ell,\C),
\end{eqnarray}
which is usually called the vacuum $\hat{\g}$-module.
It is well known ([FZ], [Lia], [Li2], [LL]) that
$V_{\hat{\g}}(\ell,0)$ has a natural vertex algebra structure.
It is also known ([Li2], [LL], cf. [FZ]) that
a module for $V_{\hat{\g}}(\ell,0)$ (as a vertex algebra)
exactly amounts to a restricted $\hat{\g}$-module of level $\ell$.

For the rest of this section we assume that
$\hat{\g}$ is a standard affine Lie algebra 
(with $\g$ a finite-dimensional simple Lie algebra
and with $\<\cdot,\cdot\>$ the normalized Killing form).
For any complex number $\ell$ not the negative dual Coxeter number of
$\g$, $V_{\hat{\g}}(\ell,0)$ equipped with a canonical conformal
vector is a vertex operator algebra (cf. [FZ]). 
For any restricted $\hat{\g}$-module $W$ of level $\ell$,
$W$ is naturally a module for  $V_{\hat{\g}}(\ell,0)$ viewed as a
vertex algebra, then $W$ is naturally a $\tilde{\g}$-module
with ${\bf d}$ acting as $\alpha-L(0)$ where $\alpha$ is any complex
number (cf. [LL]).
Denote by $L_{\hat{\g}}(\ell,0)$ the simple quotient vertex operator algebra
of $V_{\hat{\g}}(\ell,0)$.
If $\ell$ is a positive integer, it was proved ([FZ], [DL], cf. [Li2]) that
irreducible modules for $L_{\hat{\g}}(\ell,0)$ 
viewed as a vertex operator algebra
are exactly the irreducible highest weight integrable 
$\hat{\g}$-modules of level $\ell$.
Up to isomorphism, irreducible highest weight integrable 
$\hat{\g}$-modules of level $\ell$ are $L(\ell,\lambda)$
where $\lambda$ is a dominant integrable weight for $\g$
such that $\<\lambda,\theta\>\le \ell$ (see [K1]).

In [FHL], among other important results,
for a general vertex operator algebra $V$ and for $V$-modules
$W_{1},W_{2}$ and $W_{3}$, a notion of fusion
rule of type $\binom{W_{3}}{W_{1}W_{2}}$ was defined.
Furthermore, in [FZ], a conceptual method for determining fusions
was developed in terms of Zhu's algebra and this method was applied
to the case with $V=L_{\hat{\g}}(\ell,0)$.
In \cite{li-duke}, a certain analogue of the classical hom-functor
for vertex operator algebras was developed and by using this analogue
it was proved (\cite{li-duke}, Proposition 4.15) that 
the fusion rule of type $\binom{L(\ell,\nu)}{L(\ell,\lambda) L(\ell,\mu)}$
for the vertex operator algebra $L_{\hat{\g}}(\ell,0)$ 
equals the dimension of 
$\Hom _{\tilde{\g}}(L(\ell,\lambda)\otimes L(\mu)
\otimes \C[t,t^{-1}],L(\ell,\nu))$,
where $L(\mu)$ denotes the irreducible $\g$-module
of highest weight $\mu$.
Combining this with Proposition \ref{pequivalence} we immediately have:

\bp{pfusion}
Let $\ell$ be a positive integer and 
let $L(\ell,\lambda), L(\ell,\mu)$ and $L(\ell,\nu)$ be highest weight
irreducible $\hat{\g}$-modules of level $\ell$. Then 
the fusion rule of type $\binom{L(\ell,\nu)}{L(\ell,\lambda) L(\ell,\mu)}$
for the vertex operator algebra $L_{\hat{\g}}(\ell,0)$
equals the dimension of 
$\Hom_{\hat{\g}}(L(\ell,\lambda)\otimes L(\mu)(z), \overline{L(\ell,\nu)})$.
\ep

\br{r-A}
{\em It was proved in [CP2] that for a highest weight irreducible 
integrable $\tilde{\g}$-module $W$ and a finite-dimensional
irreducible $\g$-module $U$, $W\otimes U\otimes \C[t,t^{-1}]$ 
is an irreducible $\tilde{\g}$-module 
if $W$ and $U$ satisfy certain conditions.
In [A], the irreducibility of $\hat{\g}$-modules
$W\otimes U\otimes \C[t,t^{-1}]$ for certain nonintegrable
$\tilde{\g}$-modules $W$ was studied 
in terms of vertex operator algebra  $L_{\hat{\g}}(\ell,0)$
and fusion rules, and certain interesting results were obtained
in [A].}
\er


\begin{thebibliography}{FKRW}

\bibitem[A]{a1}
D. Adamovic, Vertex operator algebras and irreducibility of certain modules 
for affine Lie algebras, {\em Math. Research Letters} {\bf 4} (1997), 809-821.

\bibitem[B]{b}
R. E. Borcherds, Vertex algebras, Kac-Moody algebras, and the Monster,
{\it Proc. Natl. Acad. Sci. USA} {\bf 83} (1986), 3068-3071.

\bibitem[C]{ch}
V. Chari, Integrable representations of affine Lie algebras, 
{\em Invent. Math.} {\bf 85} (1986), 317-335.

\bibitem[CP1]{cp1}
V. Chari and A. N. Presslely, New unitary representations of loop groups, 
{\em Math. Ann.} {\bf 275} (1986), 87-104.

\bibitem[CP2]{cp2}
V. Chari and A. N. Presslely, A new family of irreducible, integrable modules 
for affine Lie algebras, {\em Math. Ann.} {\bf 277} (1987), 543-562.

\bibitem[CP3]{cp3}
V. Chari and A. N. Presslely, Integrable representations of twisted 
affine Lie algebras, {\em J. Algebra} {\bf 113} (1988), 438-464.

\bibitem[Di]{di}
J. Dixmier, {\em Enveloping Algebras}, 
(The 1996 Printing of the 1977 English Translation), 
Graduate Studies in Mathematics, Vol. II,
American Mathematical Society, 1996.

\bibitem[DL]{dl}
C. Dong and J. Lepowsky, {\it Generalized Vertex Algebras and
Relative Vertex Operators}, Progress in Math., {\bf Vol.} 112,
Birkh\"auser, Boston, 1993.

\bibitem[DLM]{dlm}
C. Dong, H.-S. Li and G. Mason, 
Regularity of rational vertex operator algebras,
{\em Adv. Math.} {\bf 132} (1997), 148-166.

\bibitem[FHL]{fhl}
I. Frenkel, Y.-Z. Huang and J. Lepowsky, On axiomatic approaches to
vertex operator algebras and modules, Memoirs Amer. Math.
Soc. {\bf 104}, 1993.

\bibitem[FLM]{flm}
I. Frenkel, J. Lepowsky and A. Meurman, {\it Vertex Operator Algebras
and the Monster}, Pure and Appl. Math., {\bf Vol. 134}, Academic Press,
Boston, 1988.

\bibitem[FZ]{fz}
I. Frenkel and Y.-C. Zhu, Vertex operator algebras associated to
representations of affine and
Virasoro algebras, {\it Duke Math. J.} {\bf 66} (1992), 123-168.

\bibitem[HL]{hl}
Y.-Z. Huang and J. Lepowsky, A theory of tensor product for
module category of a vertex operator algebra, I, {\em Selecta Mathematica},
{\bf 1} (1995), 699-756.
 
\bibitem[H]{hum}
J. Humphreys, {\em Introduction to Lie Algebras and Representation
Theory,} Springer-Verlag, New York, 1972.

\bibitem[K1]{kacbook}
V. G. Kac, Infinite Dimensional Lie Algebras,
3rd edition, Cambridge University Press, 1990.

\bibitem[K2]{kacpaper}
V. G. Kac, Constructing groups associated to 
infinite-dimensional Lie algebras, in:
{\em Infinite-dimensional groups with applications,} MSRI publ. 4, 
Springer-Verlag, 1985, 167-216.

\bibitem[LL]{ll}
J. Lepowsky and H.-S. Li,
{\em Introduction to Vertex Operator Algebras and Their Representations},
Birkh\"auser, Boston, to appear.

\bibitem[Li1]{li1}
H.-S. Li, Representation theory and tensor product theory for vertex
operator algebras, Ph.D. thesis, Rutgers University, 1994.

\bibitem[Li2]{li-local}
H.-S. Li, Local systems of vertex operators, vertex superalgebras and
modules, {\em J. Pure Appl. Alg.} {\bf 109} (1996), 143-195; 
hep-th/9406185.

\bibitem[Li3]{li-duke}
H.-S. Li, An analogue of the hom functor and a generalized nuclear democracy theorem,
{\em Duke Math. J.} {\bf 93} (1998), 73-114.

\bibitem[Li4]{li-simple}
H.-S. Li, Simple vertex operator algebras are nondegenerate,
{\em J. Algebra}, to appear.

\bibitem[Lia]{lia}
B. Lian, On the classification of simple vertex operator algebras, 
{\em Commun. Math. Phys.} {\bf 163} (1994), 307-357.

\end{thebibliography}
\end{document}